\documentclass[11pt,a4paper]{article}

\usepackage[american]{babel}
\usepackage[utf8]{inputenc}
\usepackage[T1]{fontenc}

\usepackage[a4paper,top=2.5cm,bottom=2.5cm,left=2.5cm,right=2.5cm,marginparwidth=2cm]{geometry}

\usepackage{amsmath}
\usepackage{graphicx}
\usepackage[colorinlistoftodos]{todonotes}
\usepackage[colorlinks=true, allcolors=black]{hyperref}
\usepackage{lipsum} 
\usepackage{floatrow} 
\usepackage[hypcap=false]
{caption} 
\usepackage{subcaption}

\usepackage{listings}
\usepackage{url}
\usepackage{courier}
\usepackage{wrapfig}
\usepackage{color,soul}
\definecolor{codegreen}{rgb}{0,0.6,0}
\usepackage{chngpage}
\usepackage{xcolor}
\usepackage{bm}
\usepackage{amsmath,amssymb}
\usepackage{verbatim}
\newlength{\mycolwidth}
\usepackage{array}
\newcolumntype{Z}{>{$}p{\mycolwidth}<{$}}
\usepackage{algorithmicx}
\usepackage{algorithm} 
\usepackage{refcount}
\usepackage{algpseudocode} 
\usepackage{placeins} 
\MakeRobust{\Call}

\usepackage{amsfonts}
\usepackage{graphicx}
\usepackage{epstopdf}
\usepackage{amsmath} 
\usepackage{booktabs}
\usepackage{comment}

\usepackage{stmaryrd}

\DeclareMathAccent{\svec}{\mathord}{letters}{126}
\newcommand\stvec[1]{\mathbf #1}				
				
\newcommand\ssvec[1]{\svec{\stvec{#1}}}

\usepackage{authblk}
\usepackage[labelformat=simple]{subcaption}
\usepackage{cite}

\usepackage{lineno}

\title{A comparative study of explicit and implicit Large Eddy Simulations using a high-order discontinuous Galerkin solver: application to a Formula 1 front wing}
\author[,1,2]{Gerasimos Ntoukas \thanks{Corresponding Author\newline \hspace*{1.55em} E-mail address: gerasimos.ntoukas@upm.es (G. Ntoukas)}}
\author[2,3]{Gonzalo Rubio}
\author[2]{Oscar Marino}
\author[4]{Alexandra Liosi}
\author[4]{Francesco Bottone}
\author[4]{Julien Hoessler}
\author[2,3]{Esteban Ferrer}
\affil[1]{Cadence Design Systems, Chaussée de la Hulpe 187, B-1170 Brussels, Belgium}
\affil[2]{ETSIAE-UPM - School of Aeronautics, Universidad Politécnica de Madrid, Plaza Cardenal Cisneros 3, E-28040 Madrid, Spain}
\affil[3]{Center for Computational Simulation, Universidad Politécnica de Madrid, Campus de Montegancedo, Boadilla del Monte, 28660 Madrid, Spain}
\affil[4]{CFD Methodology Group, McLaren Racing, Woking, United Kingdom}
\date{January 2024}

\begin{document}

\maketitle
\thispagestyle{empty}

\begin{abstract}

This paper explores two Large Eddy Simulation (LES) approaches within the framework of the high-order discontinuous Galerkin solver, Horses3D. The investigation focuses on an Inverted Multi-element Wing in Ground Effect (i.e. 2.5D Imperial Front Wing section) representing a Formula 1 front wing, and compares the strengths and limitations of the two LES methods. 

The explicit LES formulation relies on the Vreman model, that adapts to laminar, transitional and turbulent regimes. The numerical formulation uses nodal basis functions and Gauss points. The implicit LES formulation, does not require explicit turbulence modeling but relies in the discretization scheme. We use the Kennedy-Gruber entropy stable formulation to enhance stability in under resolved simulations, since we recover the continuous properties such as entropy conservation at a discrete level. This formulation employs Gauss-Lobatto points, which downgrades the accuracy of integration but allows for larger time steps in explicit time integration. 

We compare our results to Nektar++ \cite{slaughter2023large} showing that both LES techniques provide results that agree well with the reference values. The implicit LES shows to better capture transition and allows for larger time steps at a similar cost per iteration. We conclude that this implicit LES formulation is very attractive for complex simulations.

\end{abstract}
\mbox{}
\vfill
\newpage

\section{Introduction}
Turbulent flows play a pivotal role in engineering applications, and in particular in car aerodynamics. Understanding and accurately predicting turbulent phenomena remains a challenge in computational fluid dynamics. Large Eddy Simulation (LES) is a powerful tool for capturing the essential features of turbulent flows by resolving the larger energy-containing eddies while modeling the smaller, unresolved scales. In this context, this scientific paper delves into the comprehensive investigation of two distinct approaches for LES: Explicit LES (eLES) utilizing the Vreman subgrid-scale model and Implicit LES (iLES) incorporating the Kennedy-Gruber entropy stable formulation to enhance stability. Both methodologies are implemented within the high-order discontinuous Galerkin (DG) solver Horses3D \cite{HORSES3D}, which provides a framework for the accurate representation of turbulent flows. Horses3D is an open source Fortran 2003 object-oriented solver. In this work, we consider the compressible Navier-Stokes version using a high-order DG spatial discretizations and explicit time marching.  

Horses3D includes well-established explicit LES models, including the classic Smagorinsky model \cite{smagorinsky1963general,lilly1965computational}, the Wale model \cite{nicoud1999subgrid}, and the Vreman model \cite{vreman2004eddy}. We here focus on the Vreman subgrid-scale model (using nodal Gauss points) for its ability to provide an explicit representation of the subgrid-scale stresses, while requiring low computational cost. The explicit nature of the method allows for a direct integration of the subgrid-scale model, offering computational efficiency and ease of implementation. Our investigation into eLES with the Vreman model aims to explore its performance in capturing turbulence features for a Formula 1 front wing, shedding light on its strengths and limitations in comparison to an implicit LES formulation.

Implicit LES presents an alternative paradigm by embedding the turbulence modeling directly into the discretization scheme, eliminating the need for explicit subgrid-scale modeling \cite{buscariolo2021computational}. Numerous studies have been dedicated to developing entropy- and energy-stable schemes for the discontinuous Galerkin method, as evidenced by the extensive literature \cite{FISHER2013518,doi:10.1137/130928650,10.1016/j.jcp.2016.09.013,10.1007/s10915-018-0702-1,MANZANERO2020109241,CHEN2017427,WINTERS20181,doi:10.1137/120890144}, including comparisons between explicit and implicit LES formulations in the context of high-order methods (e.g., \cite{FERRER2019104239,KOU2023112399}, but to the author's knowledge not for the complex case of a simplified Formula 1 front wing.
 Comprehensive reviews of energy/entropy stable schemes tailored for discontinuous Galerkin schemes can be found in \cite{10.1016/j.jcp.2016.09.013,winters2021construction,CSIAM-AM-1-1}. The Kennedy-Gruber \cite{kennedy2008reduced} entropy stable formulation, renowned for its ability to maintain stability and accuracy, is employed in the implicit context. Horses3D has the capability to employ Gauss-Lobatto (GL) quadrature points (instead of Gauss points, as in the Vreman model) for integral approximation. This choice, made possible by its Summation-By-Parts Simultaneous-Approximation Term (SBP-SAT) property \cite{FISHER2013518,doi:10.1137/130932193}, facilitates the recovery of continuous properties such as energy and entropy conservation at a discrete level. This feature significantly enhances simulation robustness, particularly in scenarios characterized by high numerical aliasing, such as under-resolved turbulence \cite{DBLP:journals/jscic/ManzaneroRFVK18}, albeit at a higher computational cost.
 The iLES approach provides inherent advantages in terms of numerical stability and alleviates the need for specifying subgrid-scale model parameters, making it an attractive option for simulating turbulent flows with complex physics, such as a Formula 1 front wing. 

Horses3D serves as the common numerical platform for both explicit and implicit LES methodologies in our study.  The DG framework leverages polynomial approximations on element-wise subdomains, offering superior accuracy and flexibility in handling complex geometries. DG exhibits high-order accuracy and demonstrates favorable dispersion and dissipation characteristics. Its compact stencil, which requires only interface neighbor data, contributes to excellent parallel scaling. Additionally, the DG scheme is renowned for its proficiency in handling unstructured and non-conforming grids. The high-order accuracy of the DG method allows for more faithful representation of turbulent structures and a reduction in numerical dissipation, crucial for capturing the intricacies of turbulent flows in Formula 1 configurations. In this work, we set the polynomial order to 4 (5th order accuracy) and simulate the Formula 1 front wing with the two LES models. We compare our results to Nektar++ \cite{slaughter2023large} and show that both LES techniques provide results that agree very well with the reference values.

%

In the fast-paced world of Formula 1 racing, the quest for aerodynamic excellence is unceasing. The front wing of a Formula 1 car, with its intricate design and pivotal role in managing aerodynamic forces and flow, remains a focal point for continuous innovation and optimization. Computational simulations have become invaluable tools in the pursuit of enhanced performance, providing insights into the complex flow phenomena surrounding the front wing.
In this era of stringent regulations and the perpetual need for speed, the exploration of advanced simulation methods has become imperative for Formula 1 teams. This paper showcases the potential of the DG approaches, shedding light on its application in the context of front wing aerodynamics and its advantages in terms of accuracy and computational efficiency. Furthermore, the incorporation of entropy stable schemes ensures robust and stable simulations, making this methodology a valuable addition to the repertoire of Formula 1 aerodynamic engineers.

The following sections will provide a comprehensive overview of the methodology, computational setup, and results obtained using this innovative simulation approach. The pursuit of aerodynamic excellence in Formula 1 is a relentless endeavor, and the research presented in this paper promises to contribute significantly to the ongoing quest for performance gains through advanced computational modeling.

\section{High order discontinuous Galerkin solver: Horses3D}
The simulations have been generated using the high-order spectral element CFD solver Horses3D \cite{HORSES3D}. Horses3D is a 3D parallel framework developed at ETSIAE–UPM (the School of Aeronautics of the Polytechnic University of Madrid).
This framework uses the high-order discontinuous Galerkin spectral element method (DGSEM) and is written in modern Fortran 2003. It targets simulations of fluid-flow phenomena such as those governed by compressible Navier-Stokes equations, and it supports curvilinear hexahedral meshes of arbitrary order. All cases are run with a low storage Runge-Kutta 3 time advancement.

\subsection{Explicit Large Eddy Simulations}
The explicit LES methodology employs spatial filtering to resolve large structures while focusing the modeling efforts on small unresolved turbulent structures, such as small eddies, exhibiting isotropic behavior. Our implementations encompass well-known LES models, including those proposed by Smagorinsky \cite{smagorinsky1963general,lilly1965computational}, Wale \cite{nicoud1999subgrid}, and Vreman \cite{vreman2004eddy}. In this work, we select Vreman, which has been adapted to the high-order discretisation and computes turbulent viscosity for each high-order nodal point within each element in the mesh. We use the standard constant $c=0.07$ and a modified length scale defined as $\Delta=\frac{V^{1/3}}{N+1}$ where $V$ is the volume of the element and $N$ is the polynomial order of the approximation. The detailed formulation is given in Appendix~\ref{sec:cNS}.
The key feature of the Vreman LES model is its dynamic character, which enhances the accuracy of the simulation, involving inhomogeneous turbulence, by adjusting the model parameters based on the local flow characteristics. In particular, the Vreman model automatically reduces in laminar, transitional, and near-wall regions allowing to capture the correct physics.
This adaptability is particularly useful for capturing a wide range of turbulent behaviors in different regions of a simulation domain. The explicit LES implementation uses Gauss points.

\subsection{Implicit  Large Eddy Simulations}
LES leverages the inherent numerical dissipation originating from the numerical scheme, such as through upwinding non-linear terms, viscous terms or finite element stabilization terms, to accommodate subgrid effects. Consequently, explicit modeling of unresolved terms is unnecessary. Implicit methods have gained popularity, particularly when paired with high-order numerical techniques. These methods exhibit dissipation and dispersion errors confined to high wave number ranges, effectively constraining numerically induced dissipation to weakly-resolved regions. The term "energy stable" is associated with formulations that preserve fundamental physical properties, such as conservation of energy, and prevent the unbounded growth of numerical errors (aliasing) that can lead to instability. These formulations achieve stability by controlling the growth of certain mathematical quantities related to the simulation, often expressed as the "energy" of the system.
When simulations are conducted on coarse meshes, the cell or element sizes are larger, resulting in a less detailed representation of the physical system being simulated. This coarser discretization can lead to aliasing errors and divergence. However, an energy-stable formulation is designed to mitigate these issues by inherently damping and controlling numerical instabilities, even on coarser grids.
By maintaining stability on coarser meshes, energy-stable formulations offer enhanced robustness.

Here, we enhance the stability of the eLES computations using an "energy stable formulation". This type or formulations maintain stability even when using relatively coarse or less refined computational grids. The implict LES implementation uses Gauss-Lobatto points (including edge points) which allow SBP properties and entropy stability using the two point fluxes of Kennedy-Gruber. Note that using Gauss-Lobatto points leads to a reduced accuracy (aliasing crime) \cite{kopriva2010quadrature} but has the advantage of allowing for a higher time step (less restrictive CFL).

\section{Computational Setup}

The computational setup is based on the work presented by Slaughter et al. \cite{slaughter2023large}. The Reynolds number is $Re = 2.2 \times 10^5 $ based on the mainplane chord length (c). Since we utilise a compressible formulation of the Navier-Stokes equations, we set the Mach number to $M=0.1$. This is a compromise between the incurred restriction to the timestep size $dt$ and the maximum velocity in the domain, which should remain in the region of $M<0.3$ to avoid compressibility effects. The discretisation is based on a collocated discontinuous Galerkin method on Gauss-Lobatto nodes. Unlike the work of\cite{slaughter2023large}, where a Fourier expansion is used in the spanwise direction, Gauss (for Vreman) or Gauss-Lobatto (for KG) nodes are used to approximate the solution in all directions. This method entails a higher computational cost, but permits the solution of more general geometries, where spanwise periodicity is not viable.
The polynomial order has been set in all directions to $N=4$. Using the mesh described in Section~\ref{subsec:mesh} the total number of degrees of freedom is 9,678,750. As in \cite{slaughter2023large}, the solution has been allowed to develop for a total of 40 convection time units (CTU) which is defined as $CTU = tU/c$. The collection of statistics is performed during an additional 8 CTU.

\subsection{Domain}

\begin{figure}[H]
    \centering
    \includegraphics[width = 1.0\textwidth]{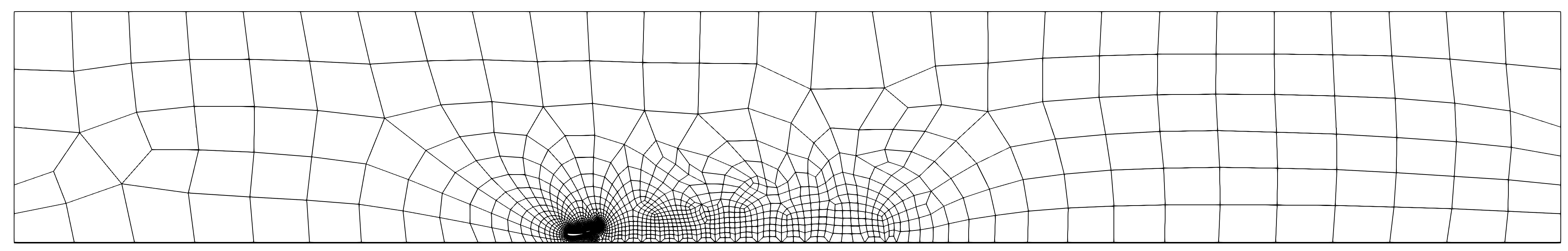}
    \caption{Overview of the domain and the mesh used for the analysis using Horses3D. The refined region is where the section of the front wing and its wake are placed. }
   \label{fig:ifw_mesh_all}
\end{figure}

The domain dimensions in the streamwise and the normal direction are augmented from the original specification of the problem. The upstream region extends 30c upstream of the leading edge of the mainplane and the outlet is positioned at 48c downstream of the trailing edge of the $2^{nd}$ flap. The height of the domain is 10c. These changes have been made to ensure that the effects of the inflow and outflow boundary conditions are diminished. Since the mesh has been designed in such way that the elements near the far-field have a significantly bigger size, as shown in Figure~\ref{fig:ifw_mesh_all}, the additional count of degrees of freedom from the extension of the domain is estimated to $\sim 1\%$. The span has been set to 0.1c which is narrower than the 0.16c used in \cite{slaughter2023large}. As shown subsequently in the results, no major differences are obtained in the obtained results. The ride height is also maintained to h/c=0.36. The tested geometry corresponds to the section of the full Imperial Front Wing at the location of $Y=250mm$. A standard inflow with $U=1.0$ and outlet are used at the ends of the domain. The top wall is also modelled as an inflow with $\Vec{u}=(1,0,0)$. The side walls of the domain use a free-slip boundary condition. In the baseline work \cite{slaughter2023large}, the authors use a periodic boundary condition, but this information was not available at the time that this work was carried out. However, this also does not seem to have a big impact as shown in the results Section~\ref{sec:Results}. The floor has been split into two sections. The section upstream of the airfoil is modeled with a free-slip wall and the rest with a moving wall using a Dirichlet boundary condition with $U=1.0$.  
The simulations are initiated under homogeneous initial conditions.

\subsection{Mesh} \label{subsec:mesh}

\begin{figure}[H]
    \centering
    \includegraphics[width = 0.85\textwidth]{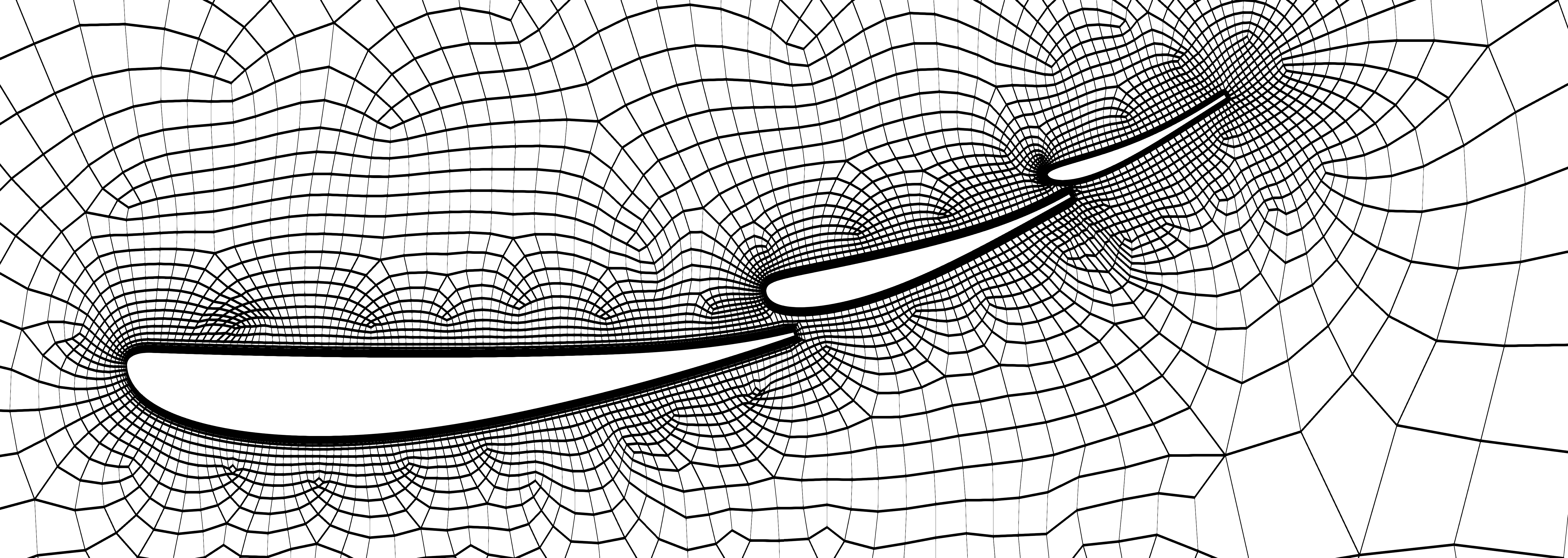}
    \caption{Close view of the mesh in the region of the three elements of the geometry. The mesh comprises of $3^{rd}$ order curvilinearl hexahedral elements.}
   \label{fig:ifw_mesh_close}
\end{figure}

An interesting aspect of this study and an important topic within the high-order methods community is that of the mesh generation. The Horses3D framework is based on DGSEM which utilizes general curvilinear hexahedral elements. Since multiple degrees of freedom are present within each element, 5 in each direction for $N=4$, the elements size is significantly bigger than those found in the Finite Volume community. This can be problematic near geometries with high-curvature. Therefore, we resort to the use of $4th$ order curvilinear elements to approximate the geometry.\\
The baseline mesh has been created using STAR-CCM+ \cite{ccm0u}. Initially a linear 2D quad-based mesh has been made which has been extruded in the spanwise direction. This mesh has been imported in the HOPR library \cite{HOPR}, an open source mesh generator which is able to produce a curved mesh based on the original linear mesh. This is a straightforward process, which did not require multiple iterations or any form of fine tuning, that produces a high-quality mesh suitable for high-order methods. An overiew of the mesh and a close-up near the surface are presented in Figures~\ref{fig:ifw_mesh_all} and \ref{fig:ifw_mesh_close}. \\
The mesh has been constructed such that a $N\geq 4$ polynomial order is used. As denoted in \cite{slaughter2023large}, this is a peculiar case due to the high $-C_{p_{max}}$ and the compression of boundary layer that is causes and the complex phenomena associated with the transition mechanism near the elements surface. Therefore, high resolution near the surface is of paramount importance to get a physical solution. The first element height is $h_{1}/c= 1.082 \times 10^{-4}$. A total of 12 number of layers is inserted in the near wall region with a growth factor of 1.2. Unlike, typical meshes for RANS applications, for WRLES the aspect ratio near the surface should be kept low. The aspect ratio in the streamwise direction is kept below 10. In the spanwise direction, we use 10 elements or 50 GL quadrature points. \\

\begin{figure}[H]
\centering
\begin{subfigure}{.37\textwidth}
  \centering
  \includegraphics[width = 1\textwidth]{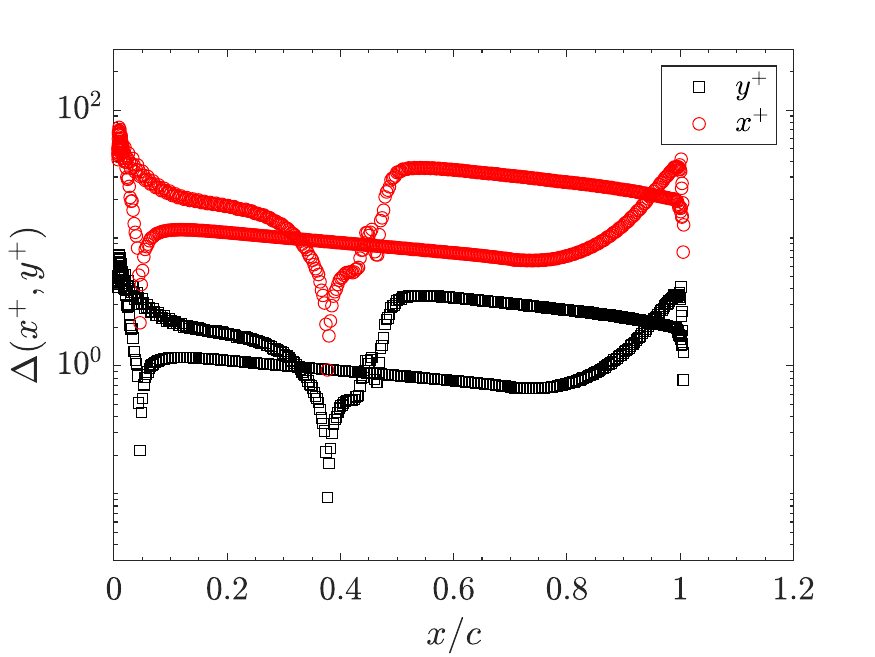}
  \caption{Mainplane}
  \label{fig:yplus_wing1}
\end{subfigure}%
\begin{subfigure}{.37\textwidth}
  \centering
  \includegraphics[width = 1\textwidth]{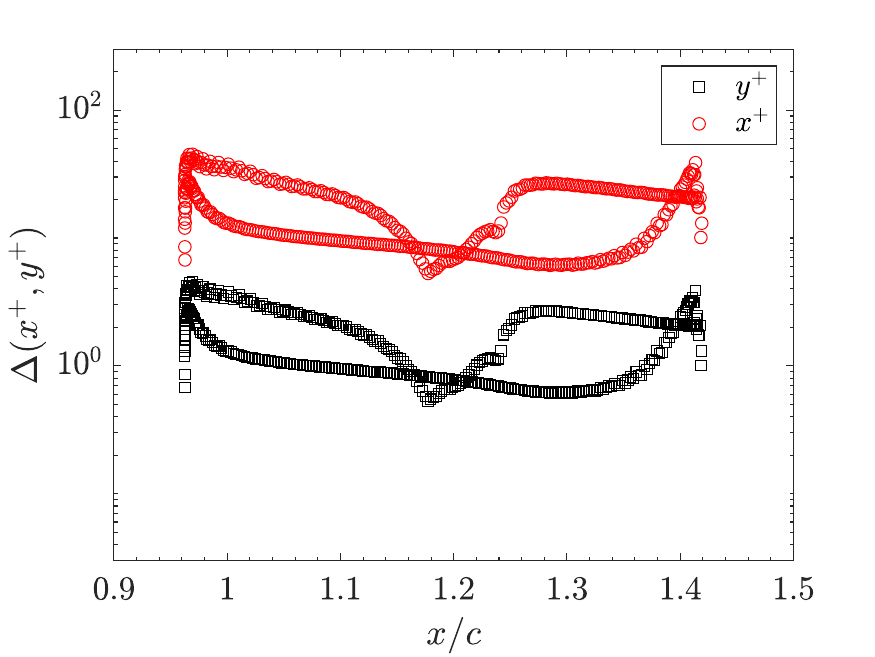}
  \caption{First flap}
  \label{fig:yplus_wing2}
\end{subfigure}
\vfill
\begin{subfigure}{.37\textwidth}
  \centering
  \includegraphics[width = 1\textwidth]{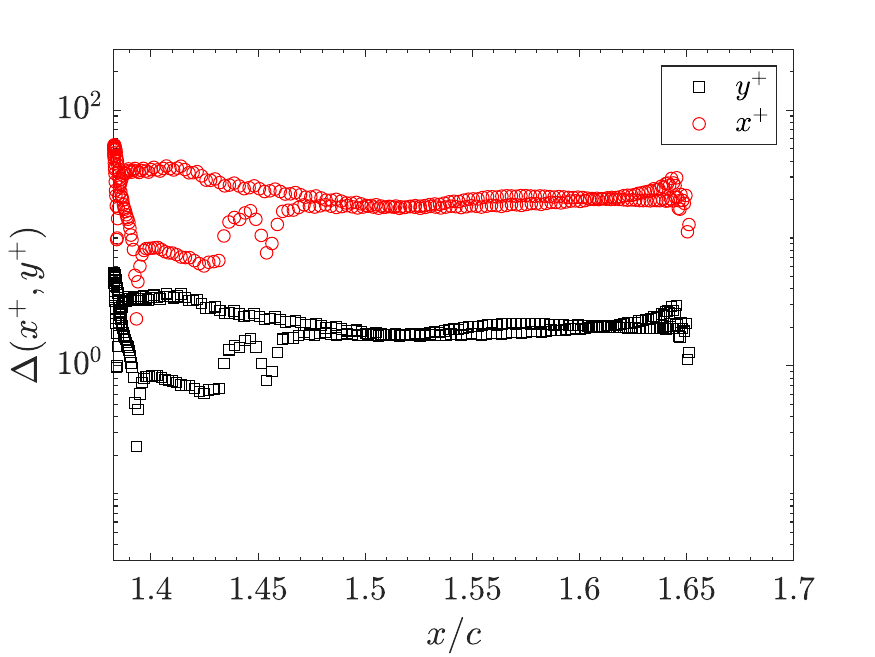}
  \caption{Second flap}
  \label{fig:yplus_wing3}
\end{subfigure}
\caption{$ X^{+}$ and $Y^{+}$ values for the iLES simulations using Horses3D for ($\mathrm{\protect\subref{fig:yplus_wing1}}$) the mainplane, ($\mathrm{\protect\subref{fig:yplus_wing2}}$) the first flap  and ($\mathrm{\protect\subref{fig:yplus_wing3}}$) the second flap.}
\label{fig:ifw_yplus_xplus}
\end{figure}

Figure~\ref{fig:ifw_yplus_xplus} presents the obtained $Y^{+}$ and $X^{+}$ values for each of the three elements. For the $y^{+}$ results presented in \ref{fig:ifw_yplus_xplus}, it is worth mentioning that we measure the distance to the second GL point. The first GL point is on the element boundary and thus we have placed two computational points in the linear sublayer and are able to accurately compute the near wall gradients. The $X^{+}$ values are kept below 110 which are the recommendations included in \cite{slaughter2023large}. During experimentation with different near wall resolutions, coarser grids lead to nonphysical solutions which either comprised of large separation regions or different transition mechanisms. Therefore, the near-wall resolution is extremely important and the guidelines should be strictly followed. 
\begin{figure}[H]
\centering
\begin{subfigure}{.37\textwidth}
  \centering
  \includegraphics[width = 1\textwidth]{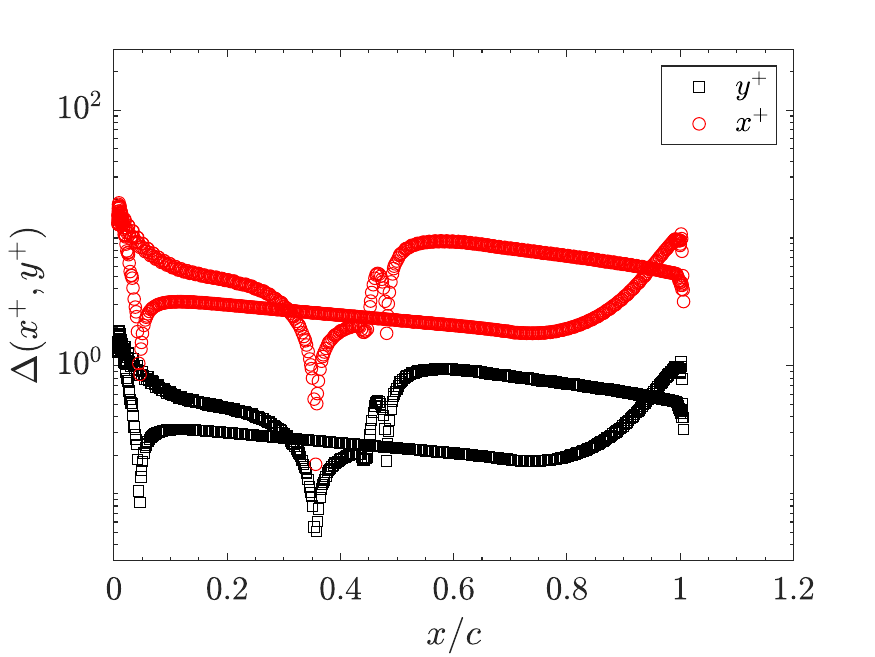}
  \caption{Mainplane}
  \label{fig:yplus_wing1_vrm}
\end{subfigure}%
\begin{subfigure}{.37\textwidth}
  \centering
  \includegraphics[width = 1\textwidth]{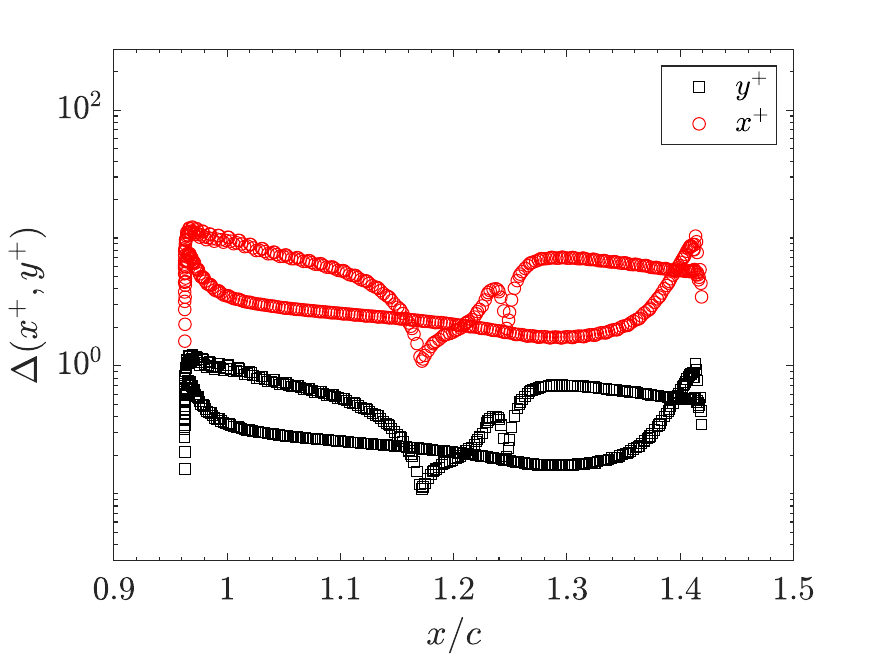}
  \caption{First flap}
  \label{fig:yplus_wing2_vrm}
\end{subfigure}
\vfill
\begin{subfigure}{.37\textwidth}
  \centering
  \includegraphics[width = 1\textwidth]{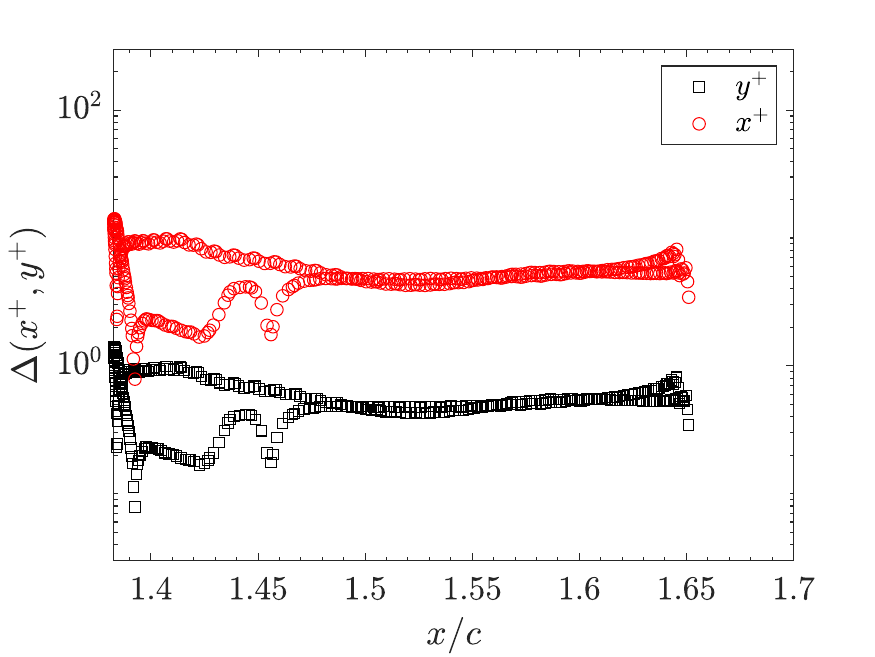}
  \caption{Second flap}
  \label{fig:yplus_wing3_vrm}
\end{subfigure}
\caption{$X^{+}$ and $Y^{+}$ values for the eLES simulations using Horses3D for ($\mathrm{\protect\subref{fig:yplus_wing1_vrm}}$) the mainplane, ($\mathrm{\protect\subref{fig:yplus_wing2_vrm}}$) the first flap and ($\mathrm{\protect\subref{fig:yplus_wing3_vrm}}$) the second flap.}
\label{fig:ifw_vrm_yplus_xplus}
\end{figure}

\section{Results}
\label{sec:Results}
In this section we perform a comparison between the explicit and implicit LES results obtained in this study using Horses3D. We also include the results presented in \cite{slaughter2023large} which have been obtained using Nektar++ \cite{cantwell2015nektar++,moxey2020nektar++}. Some small deviation is expected due to the differences in the setup, as outlined in Section~\ref{sec:averaged_data}, however, it is expected that the results from Horses3D should be in good agreement with those of Nektar++. Initially, we focus on averaged quantities and the flow characteristics near the surface of each element. This will allow to evaluate the prediction quality of the near-wall flow and especially how the transition mechanism is captured. Then, in \ref{sec:integral_values}, we examine the time-averaged integral quantities and present the results for the PSD of the lift coefficient for the geometry as well as for each element individually. The last part of the analysis in \ref{sec:wake_stats}, entails a comparison of different wake statistics quantities. In each subsection we evaluate the quality of the solution for each type of LES simulation such that we can assess the strengths and weaknesses of the numerical setup for this test case. 

\subsection{Averaged data}
\label{sec:averaged_data}
\begin{figure}
    \centering
    \begin{subfigure}{0.85\textwidth}
        \includegraphics[width=\textwidth]{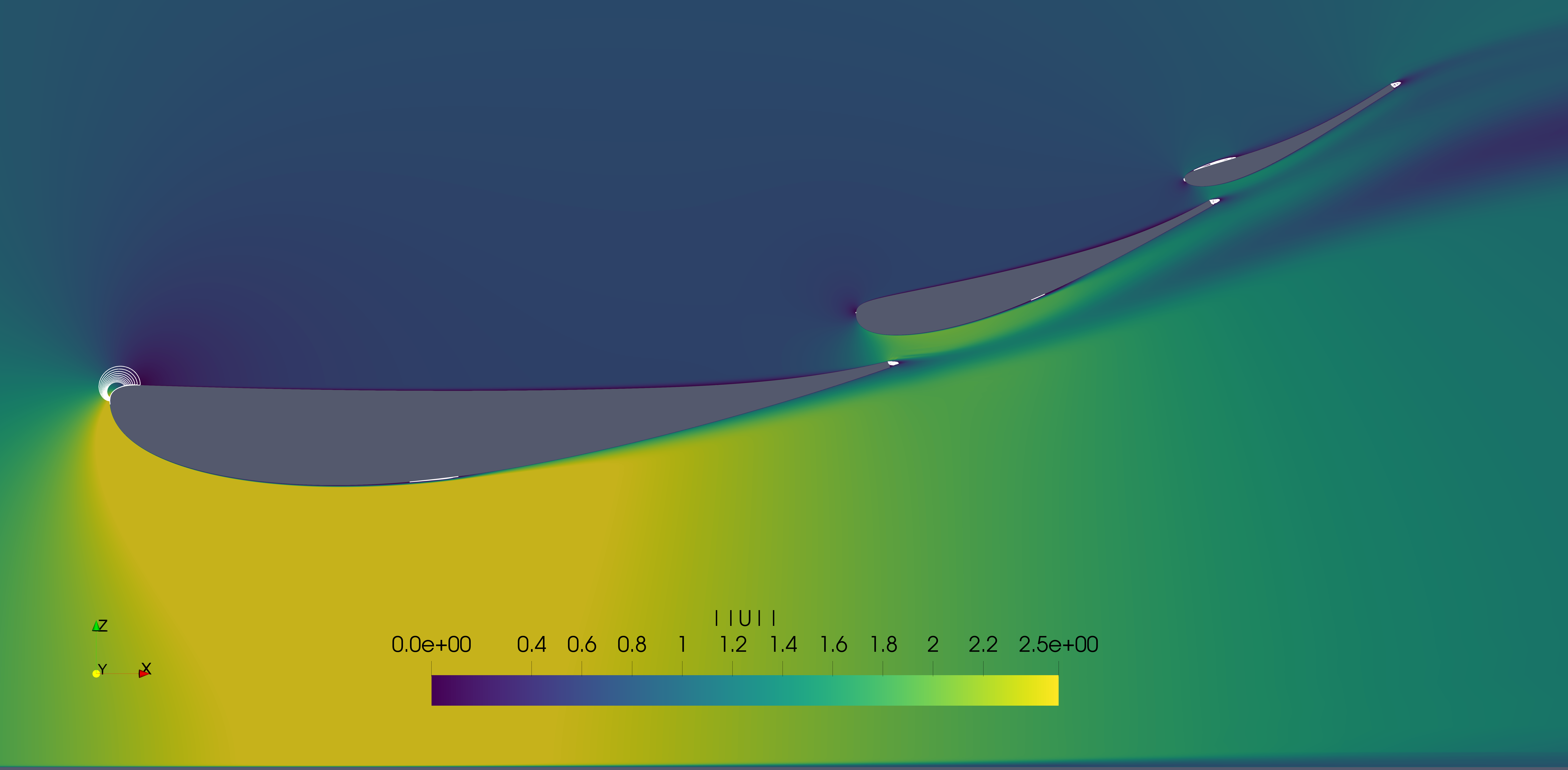}
        \caption{iLES - Horses3D}
        \label{fig:uvel_avg_Horses3D}
    \end{subfigure}
    \vfill
        \begin{subfigure}{0.85\textwidth}
        \includegraphics[width=\textwidth]{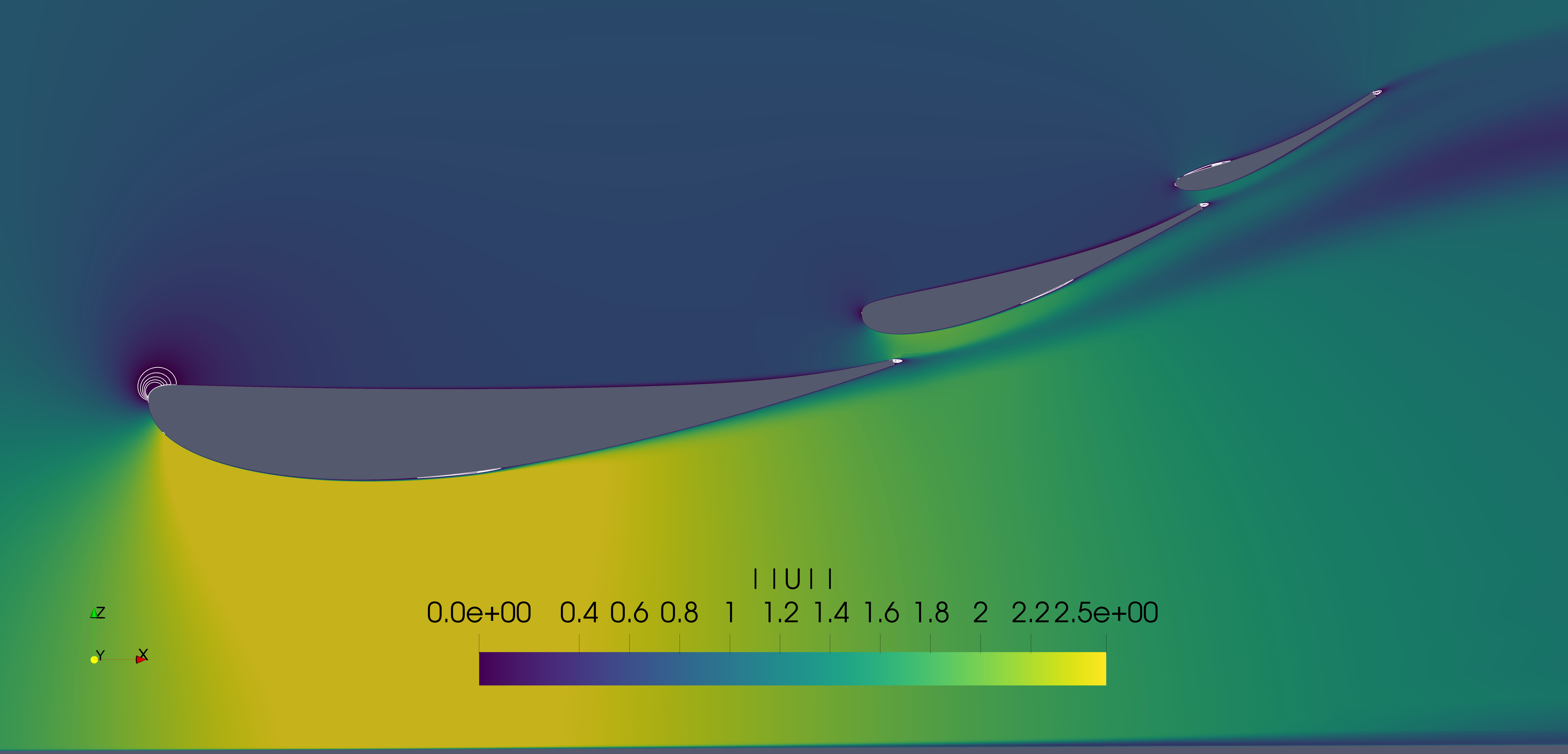}
        \caption{eLES - Horses3D}
        \label{fig:uvel_avg_vreman_Horses3D}
    \end{subfigure}
    \vfill
    \begin{subfigure}{0.8\textwidth}
        \includegraphics[width=\textwidth]{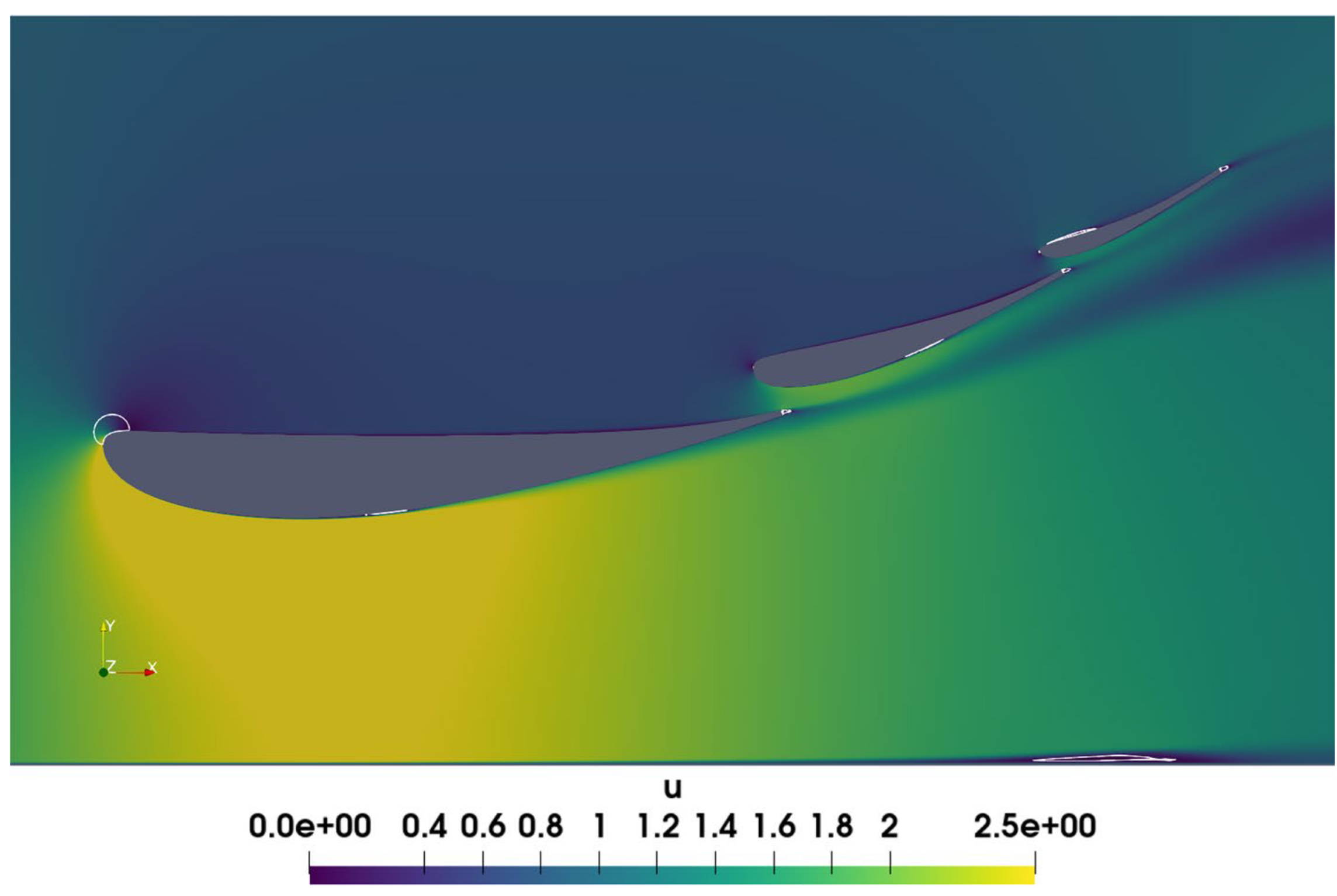}
        \caption{Nektar++ \cite{slaughter2023large}}
        \label{fig:uvel_avg_nektar}
    \end{subfigure}
    \caption{Comparison of time-averaged velocity magnitude contours for the ($\mathrm{\protect\subref{fig:uvel_avg_Horses3D}}$) iLES results, ($\mathrm{\protect\subref{fig:uvel_avg_vreman_Horses3D}}$) the eLES results and for the ($\mathrm{\protect\subref{fig:uvel_avg_nektar}}$) Nektar++ results. The white isolines indicate the regions of recirculating flow.}
    \label{fig:combined_uvel_avg}
\end{figure}

The initial point of comparison pertains to the averaged flow field, illustrated in Fig. \ref{fig:combined_uvel_avg}, where we contrast the averaged velocity flow field obtained through simulations in Horses3D and Nektar++. The overall flow topology is very similar with all the main features shown in Figure~\ref{fig:uvel_avg_nektar} being present in the iLES solution in Figure~\ref{fig:uvel_avg_Horses3D} and the eLES solution of Figure~\ref{fig:uvel_avg_vreman_Horses3D}. Namely, the presence of laminar separation bubbles (LSBs) on the suction side (bottom side) of the mainplane and of the first flap (see white contours depicting the bubble). These bubbles appear elongated and exceptionally thin, owing to the slender nature of the incoming laminar boundary layer. The thinning is attributed to the locally elevated Reynolds numbers induced by positive streamwise pressure gradient associated to the ground effect \cite{slaughter2023large}. A similar transition mechanism can be observed on the pressure side of the second flap. The LSB region captured from the eLES simulation seems to be elongated compared to that of iLES and the Nektar++ solutions. This is the case for the LSB on the mainplane and the first flap. Lastly, a recirculation region is observed on the floor for the solution from Nektar++ (bottom figure). Notably, this feature is absent in the solution from Horses3D, primarily due to variations in resolution near the bottom wall.

\begin{figure}[H]
    \centering
    \includegraphics[width = 0.8\textwidth]{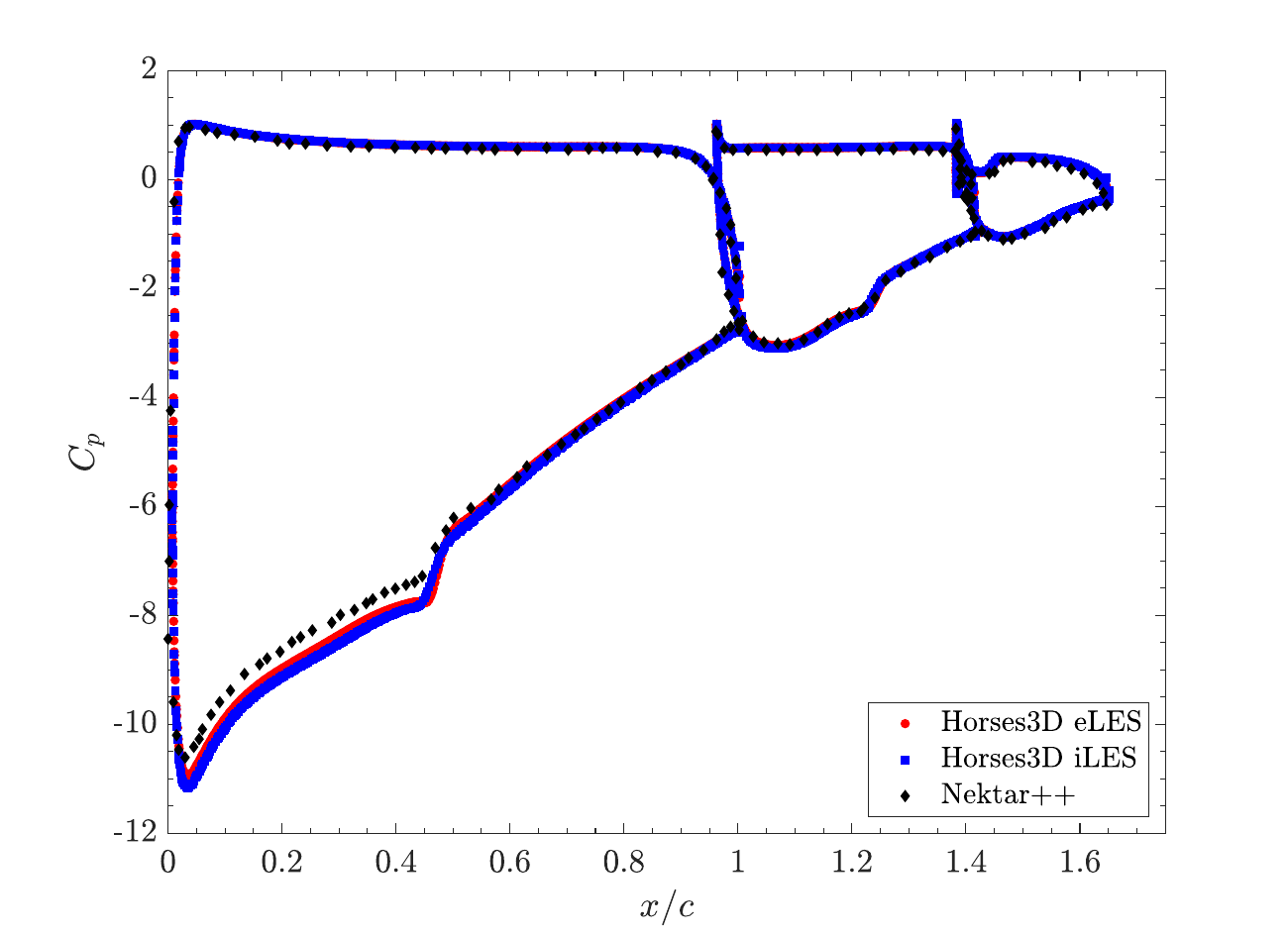}
    \caption{Average pressure coefficient $C_{p}$ comparison between Horses3D and Nektar++ for all elements of the wing section.}
   \label{fig:Cp_all}
\end{figure}

\begin{figure}[H]
    \centering
    \begin{subfigure}{0.49\textwidth}
        \includegraphics[width=\textwidth]{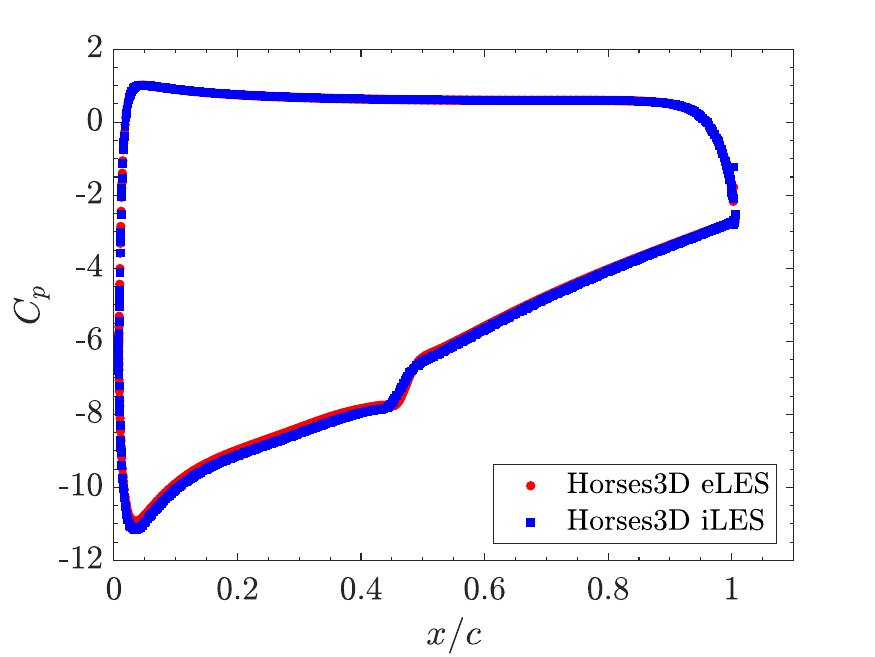}
        \caption{Mainplane}
        \label{fig:Cp_comp_iles_les}
    \end{subfigure}
    \hfill
    \begin{subfigure}{0.49\textwidth}
        \includegraphics[width=\textwidth]{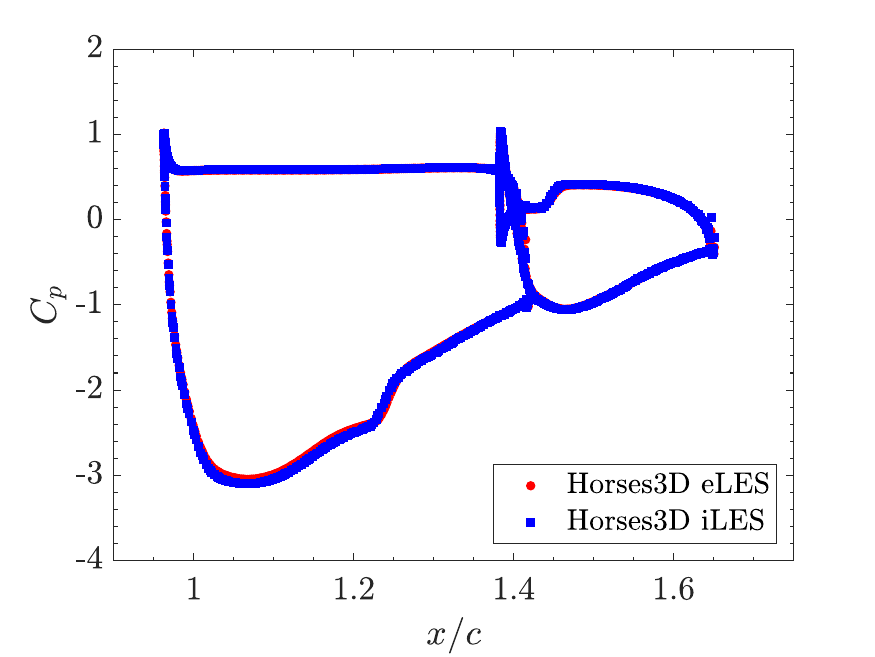}
        \caption{First and second flaps}
        \label{fig:Cp2_comp_iles_les}
    \end{subfigure}
    \caption{Comparison of average pressure coefficient $C_p$ for the iLES and eLES methodologies of Horses3D.}
    \label{fig:Cp_combined}
\end{figure}

\begin{figure}[H]
    \centering
    \includegraphics[width = 0.7\textwidth]{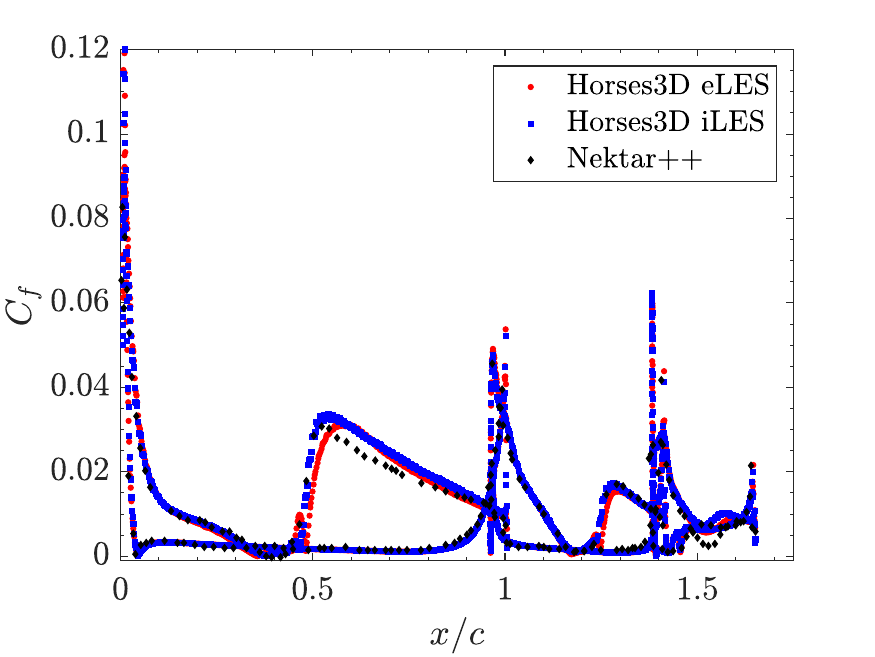}
    \caption{Average skin friction coefficient $C_{f}$ comparison between Horses3D and Nektar++ for all elements of the wing section.}
   \label{fig:Cf}
\end{figure}

\begin{figure}[H]
    \centering
    \includegraphics[width = 0.75\textwidth]{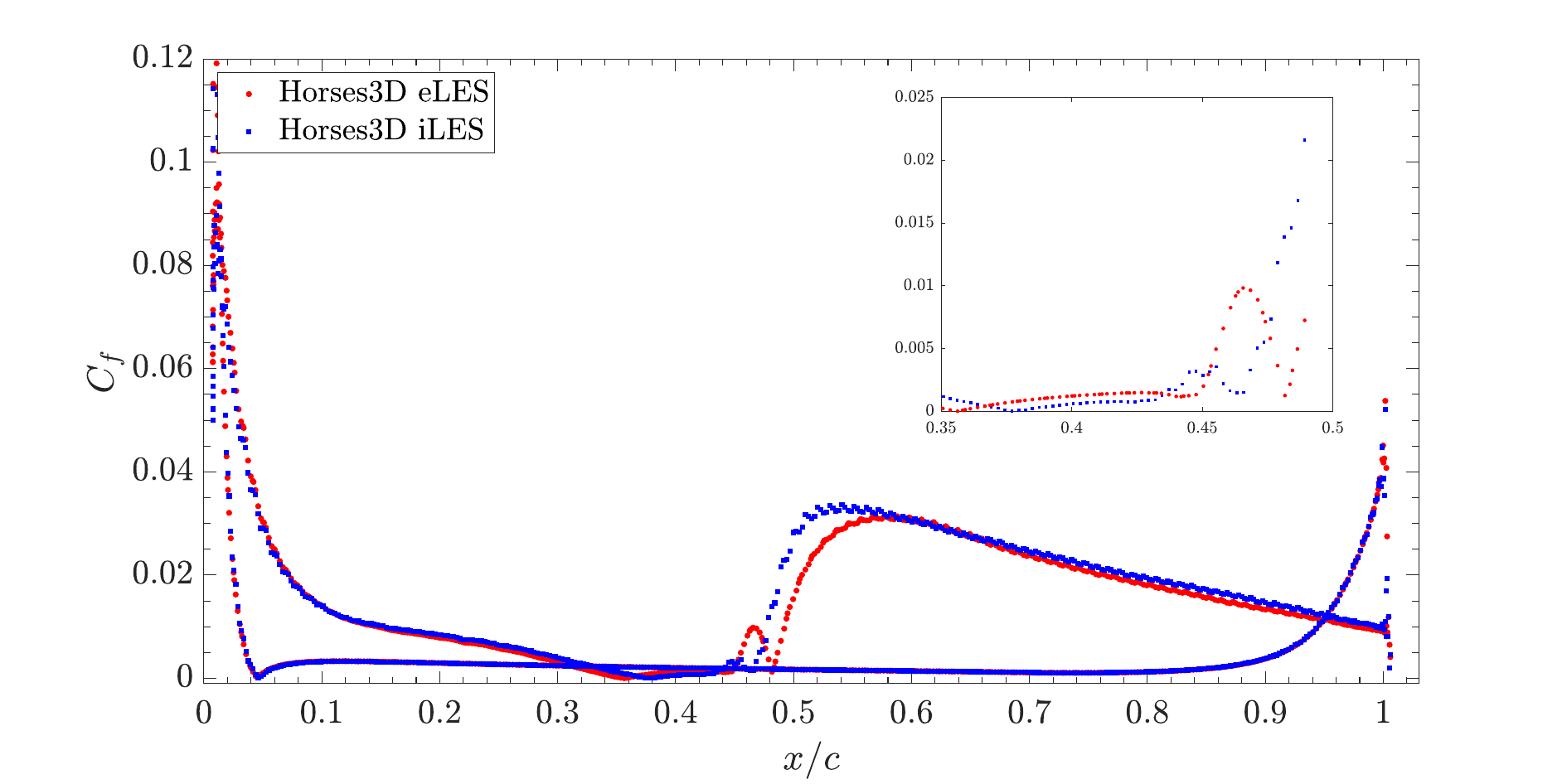}
    \caption{Average skin friction coefficient $C_{f}$ comparison between iLES and LES for the mainplane. We also present a close-up view of the $C_{f}$ around the LSB region.}
   \label{fig:Cf_comp_iles_les}
\end{figure}

\begin{figure}[H]
    \centering
    \begin{subfigure}{0.49\textwidth}
        \includegraphics[width=\textwidth]{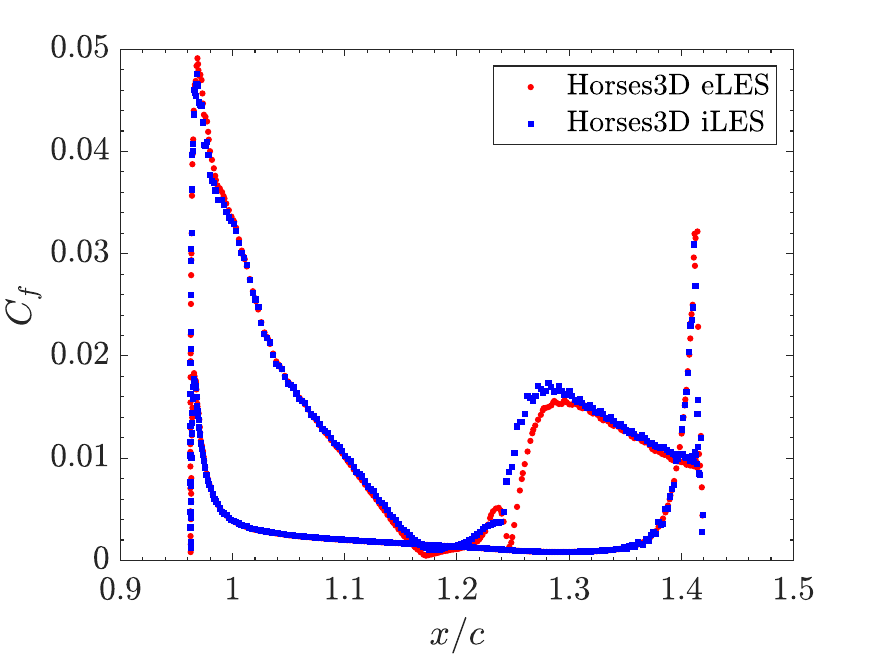}
        \caption{First flap}
        \label{fig:Cf_w2_comp_iles_les}
    \end{subfigure}
    \hfill
    \begin{subfigure}{0.49\textwidth}
        \includegraphics[width=\textwidth]{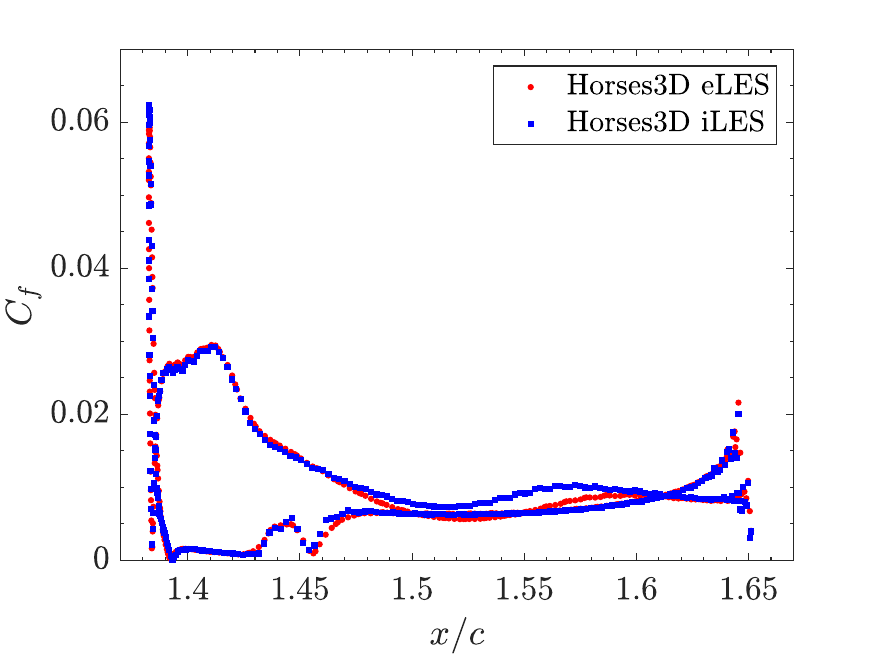}
        \caption{Second flap}
        \label{fig:Cf_w3_comp_iles_les}
    \end{subfigure}
    \caption{Average skin friction coefficient $C_{f}$ comparison between iLES and LES from Horses3D for the first and second flaps.}
    \label{fig:Cf_w2w3_comp_iles_les}
\end{figure}

\begin{figure}[H]
    \centering
    \begin{subfigure}{0.49\textwidth}
        \includegraphics[width=\textwidth]{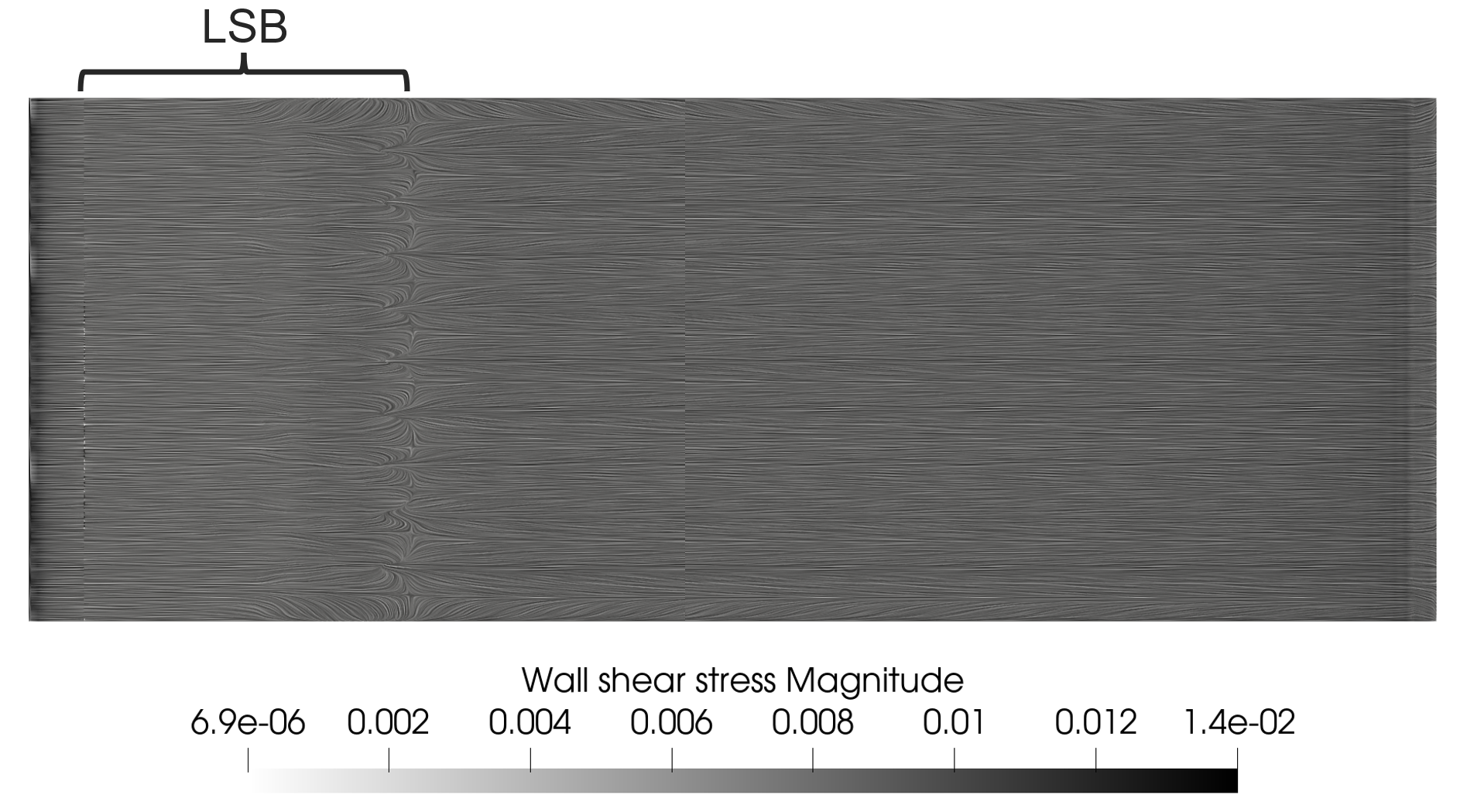}
        \caption{iLES}
        \label{fig:Cf_lic_kg_wing3}
    \end{subfigure}
    \hfill
    \begin{subfigure}{0.49\textwidth}
        \includegraphics[width=\textwidth]{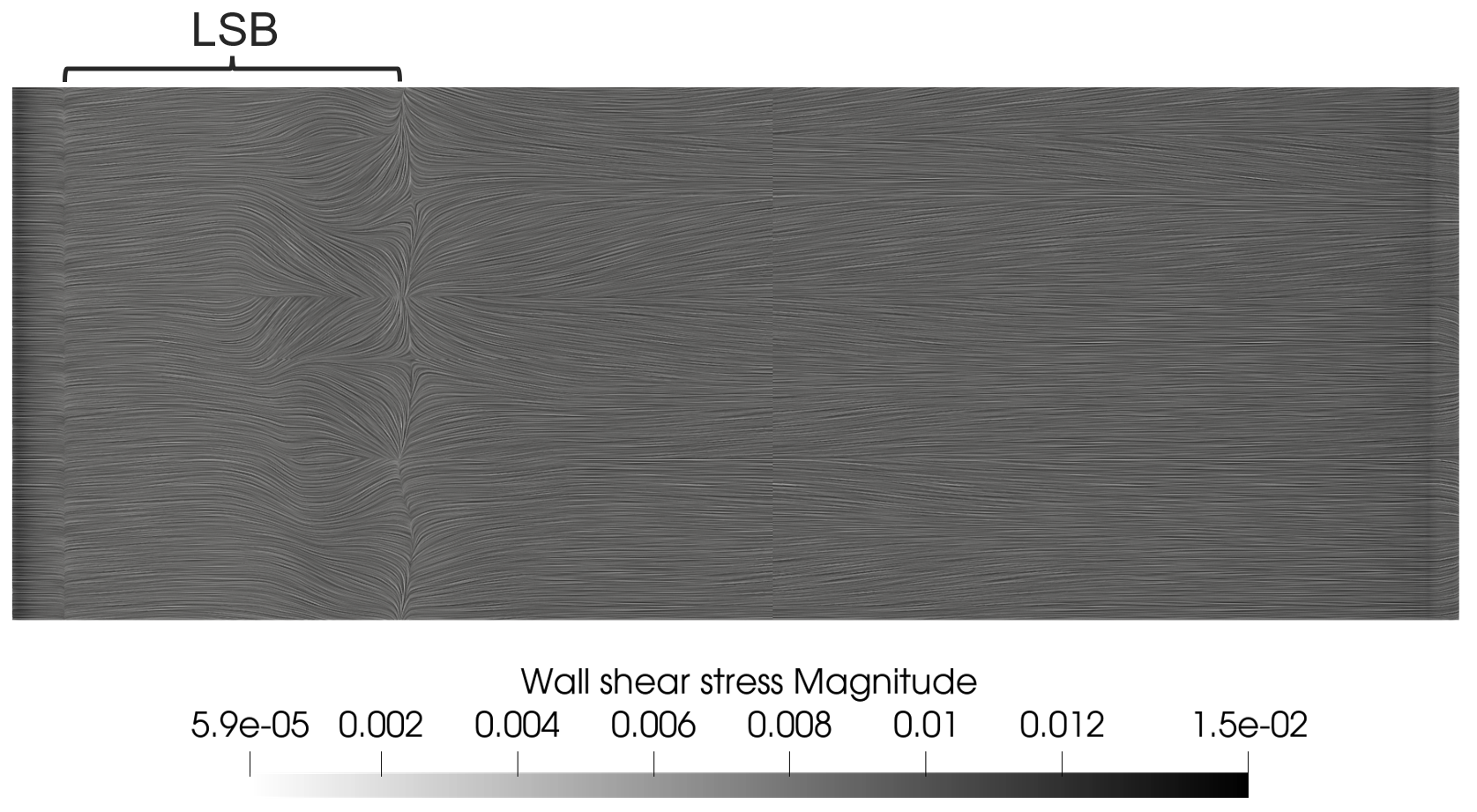}
        \caption{eLES}
        \label{fig:Cf_lic_vrem_wing3}
    \end{subfigure}
    \caption{LIC plot of the averaged wall shear stress for the second flap for the iLES and eLES results.}
    \label{fig:Cf_lic_wing3}
\end{figure}

\begin{figure}[H]
    \centering
    \includegraphics[width = 0.7\textwidth]{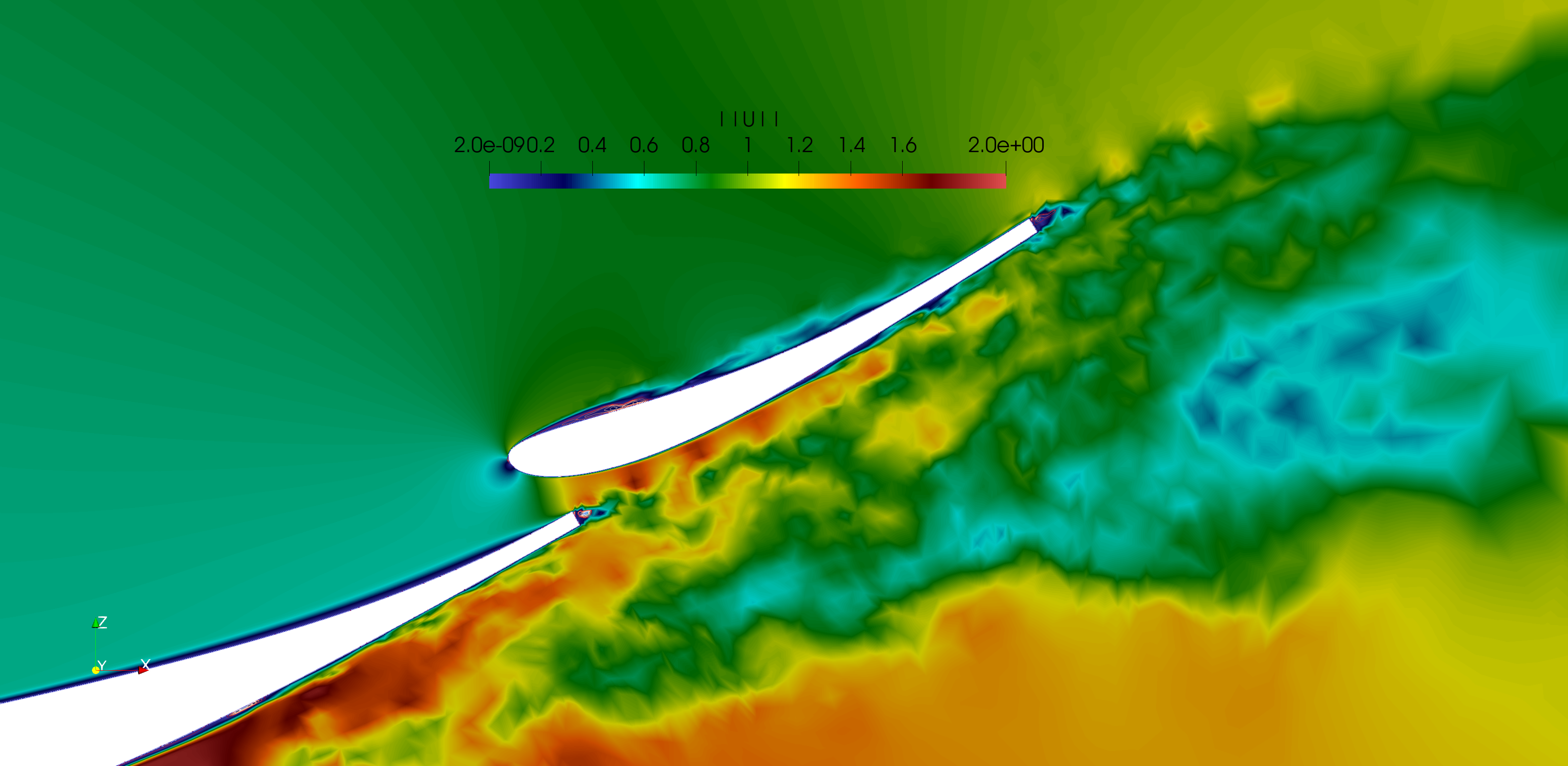}
    \caption{Instantaneous absolute velocity near the second flap. Shows the interaction of the flow on the element's surface with the wake of other elements.}
   \label{fig:uabs_inst_wing3}
\end{figure}

\begin{figure}[H]
    \centering
    \includegraphics[width = 0.8\textwidth]{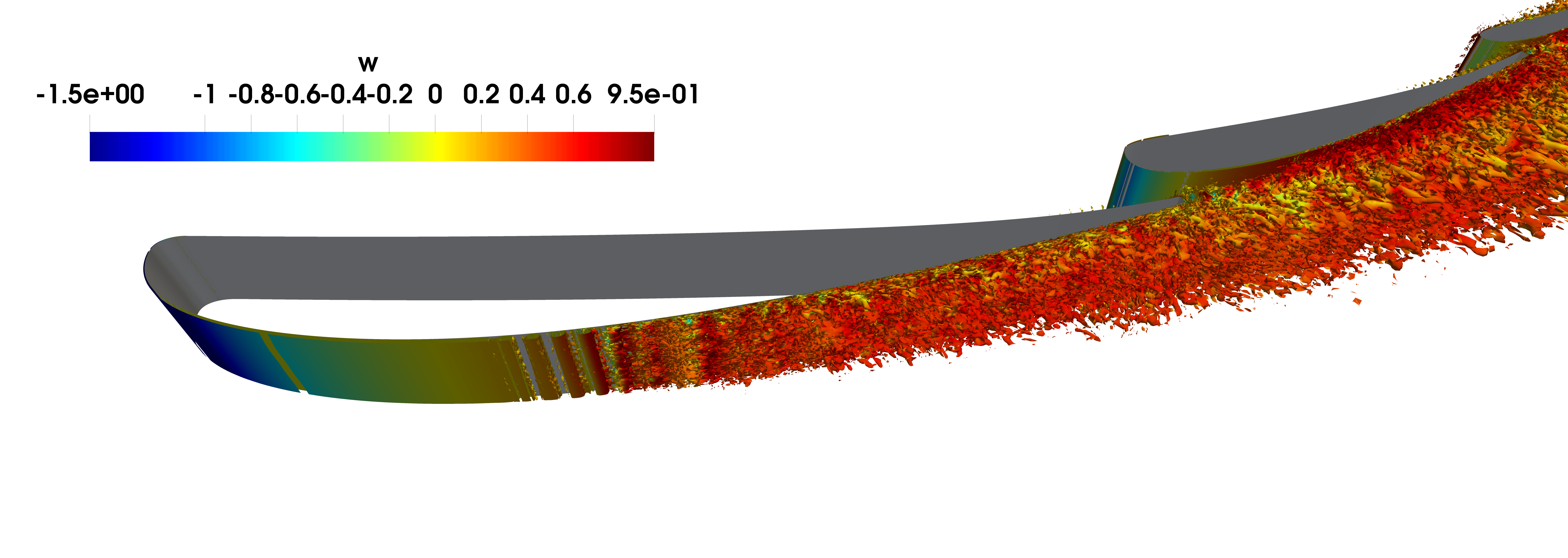}
    \caption{Instantaneous isorfaces of the Q-criterion for $Q_{crit}=100$.}
   \label{fig:Qcrit_inst_wing1}
\end{figure}

Figure \ref{fig:Cp_all} illustrates the distribution of the pressure coefficient ($C_p$) across the individual elements, defined by the equation:
\begin{equation}
    C_p = \frac{p - p_\infty}{\frac{1}{2} \rho U_\infty^2}.
\end{equation}
Upon comparison of the results obtained from Horses3D and Nektar++, the three datasets are in good agreement. The two solutions obtained from Horses3D are identical, as presented in Figure~\ref{fig:Cp_combined}, whereas some small differences are observed in the suction side of the mainplane with respect to the results from \cite{slaughter2023large}. The location of the stagnation pressure is accurately predicted. This is placed at the leading edge for the mainplane and the first flap. As noted in \cite{slaughter2023large}, the stagnation point of the second flap is not at the leading edge due to direction of the local flow with respect to that element. This shift is also present in our data. Another characteristic that we are able to reproduce is the steep adverse pressure gradient on the suction side of each element, which is followed by a monotonic increase of the $C_p$ up until the transition region. There is also a good agreement for the maximum |$C_p$| on the suction side of the mainplane and the two flaps. For the case of the mainplane, the |$C_p$| curve is slightly offset from the $C_{p_{max}}$ location up to $x=0.45$, with the solutions from Horses3D showcasing higher values. This is attributed to the use of free-slip walls on the sides of the domain and the shortest span length. For the first and second flaps we observe a very good match of the two solutions. Notably, an important aspect is the prediction of the transition region. There is a flattening and a sharp rise of the $C_p$ in the region $0.45\leq x\leq0.5$ for the mainplane and $1.2\leq x\leq1.25$ for the first flap which indicates that transition is taking place. The extent of this region can become more apparent by examining other quantities such as the skin friction coefficient or the wall shear stress.
Next we examine the average skin friction coefficient profile across the three elements, which is defined as
\begin{equation}
C_f = \frac{2 \tau_w}{\rho U_\infty^2}.
\end{equation}
The results overall seem to have a very good agreement as shown in Figure~\ref{fig:Cf}. The peaks of the $C_{f}$ near the leading edge of each element are correctly predicted and are in par with the results from Nektar++ \cite{slaughter2023large}.There is a good match of the solutions in the pressure and suction sides of the mainplane and the first flap. In the suction side, there is a rapid decrease of the $C_{f}$ and a subsequent flattening of the $C_{f}$ in the region of the LSB. These regions are very well captured from Horses3D and the results match those from Nektar++. This is followed by a sharp increase of the $C_{f}$ due to transition. Some minor differences between the eLES, iLES and Nektar++ simulations are observed in that region for the mainplane. The results from Horses3D have a slightly higher maximum $C_{f}$. For the $2^{nd}$ flap we observe a higher peak in the $C_{f}$ at the leading edge, however the region of the LSB at the pressure side is well captured. These results indicate that for this mesh, bot iLES and eLES methodologies of Horses3D are able to predict the expected flow features in the correct locations. 
In the case of the $C_{f}$, we observe some differences between the iLES and eLES in the LSB regions. To understand the differences we present the obtained $C_{f}$ for the two methods for the mainplane in Figure~\ref{fig:Cf_comp_iles_les}. The predictions for the $C_{f}$ match in the pressure side and in the portion of the suction side outside of the LSB. In the LSB region, the recirculation region extends for $0.35 \leq x/c \leq 0.48$ for the eLES and from $0.38 \leq x/c \leq 0.46$ for the iLES. A similar pattern is observed in Figure~\ref{fig:Cf_w2_comp_iles_les} for the first flap where a similar transition mechanism in the suction side exists. The extent of the LSB in the iLES results is consistent with that of \cite{slaughter2023large}, whereas the eLES solution predicts the presence of a longer LSB in both the mainplane and the first flap. For the second flap, as shown in Figure~\ref{fig:Cf_w3_comp_iles_les}, the iLES and eLES results are in agreement and the location of the LSB coincides, although in this case the LSB is created because of the orientation of the element with respect to the local flow direction.
\begin{figure}[H]
    \centering
    \includegraphics[width = 0.9\textwidth]{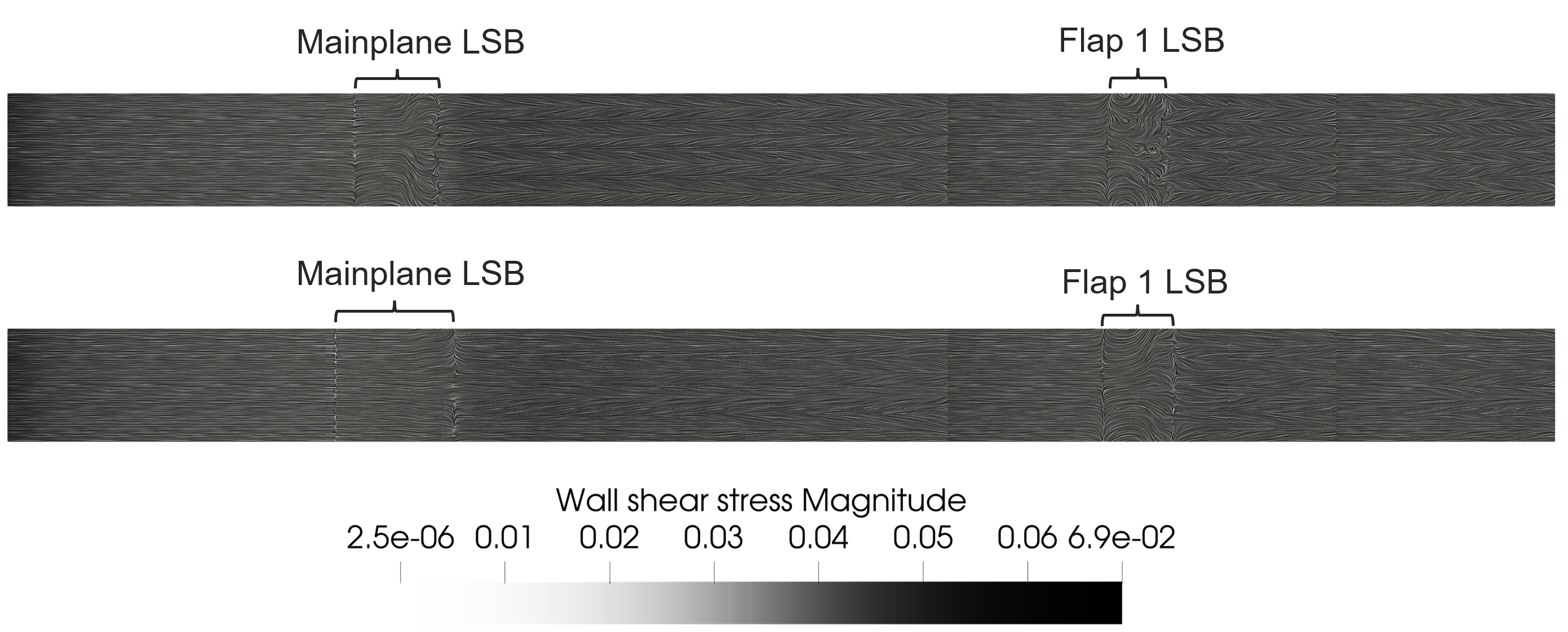}
    \caption{ LIC plot of the averaged wall shear stress. Bottom view of the three elements for the iLES results (top) and the eLES results (bottom). We have annotated the bifurcations that indicate the boundaries of the LSB for each case.}
   \label{fig:Cf_lic_all}
\end{figure}
Another very interesting way to visualise and analyse the data is a surface line integral convolution (LIC) as presented in \cite{slaughter2023large}. Figure \ref{fig:Cf_lic_all} presents the LIC of the wall shear stress on the suction side of the three elements for the iLES (top) and eLES (bottom). In the mainplane and the $1^{st}$ flap we can visually identify the location of the LSB. The difference in the length of the LSB is evident for the first two elements, with the eLES prediction being approximately $50\%$ longer than that of the iLES. For the case of the $2^{nd}$ flap no bifurcation lines are present in the suction side. All these observations match those made in \cite{slaughter2023large}. 
In Figure~\ref{fig:Cf_lic_wing3} we present the LIC for the pressure side of the $2^{nd}$ flap. In this case, the LSB is clearly visible extending from the leading edge and up to approximately a third of the chord length. The iLES and eLES results match in this case. As stated in \cite{slaughter2023large}, the high adverse pressure gradient and the interaction with the wake might be the reason why this LSB has different characteristics than those on the other elements. The effective angle of attack might also have an impact on this. A snapshot of the instantaneous absolute velocity in the region near the $2^{nd}$ flap is presented in Figure~\ref{fig:uabs_inst_wing3}, where the interaction of the $2^{nd}$ flap with the wake of the mainplane and $1^{st}$ flap can be seen.
In Figure~\ref{fig:Qcrit_inst_wing1} we present the iso-surface of the Q-criterion for $Q_{crit}=100$. After $x/c=0.5$ the flow has transitioned to fully turbulent. The structure of this mechanism can be seen on the surface of the mainplane with the formation of the LSB and the subsequent ingestion intp the subsequent burst and injection into the high momentum flow outside the boundary layer.  
\subsection{Integral values}
\label{sec:integral_values}
\begin{figure}[H]
    \centering
    \begin{subfigure}{0.49\textwidth}
        \includegraphics[width=\textwidth]{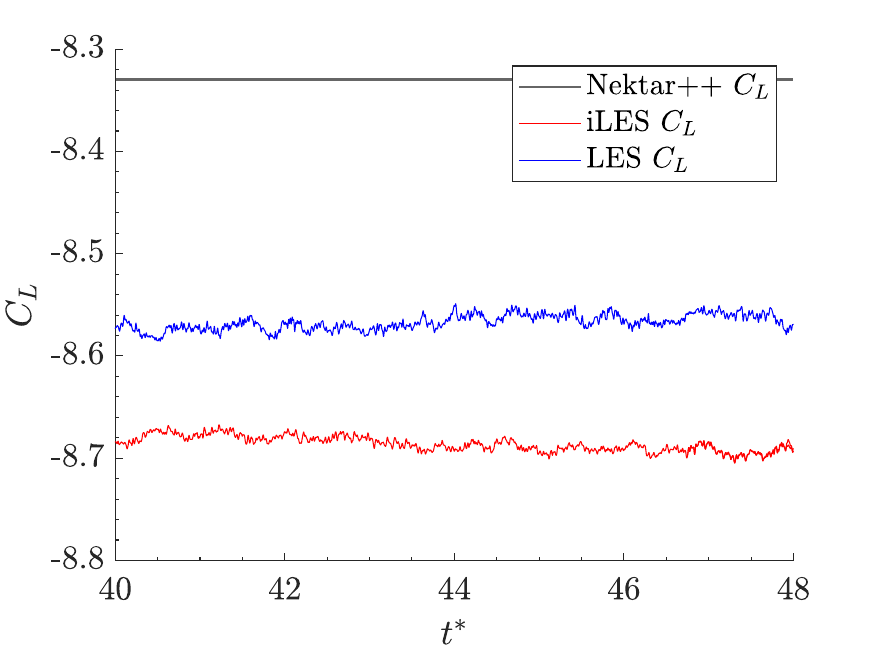}
        \caption{Lift coefficient.}
        \label{fig:lift_coeff}
    \end{subfigure}
    \hfill
    \begin{subfigure}{0.49\textwidth}
        \includegraphics[width=\textwidth]{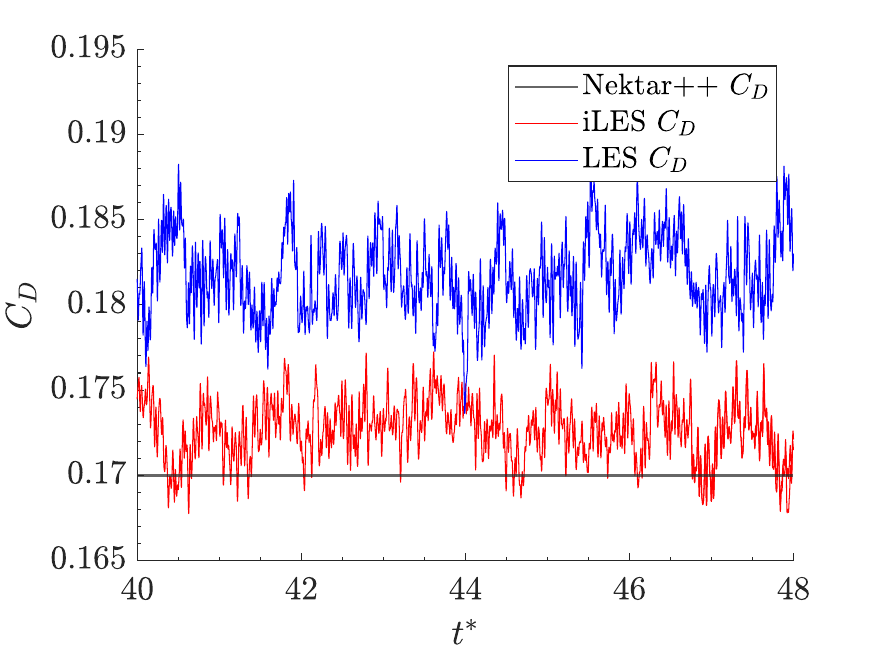}
        \caption{Drag coefficient.}
        \label{fig:drag_coeff}
    \end{subfigure}
    \caption{integral coefficients.}
    \label{fig:integral_forces}
\end{figure}

\begin{figure}[H]
    \centering
    \includegraphics[width = 0.95\textwidth]{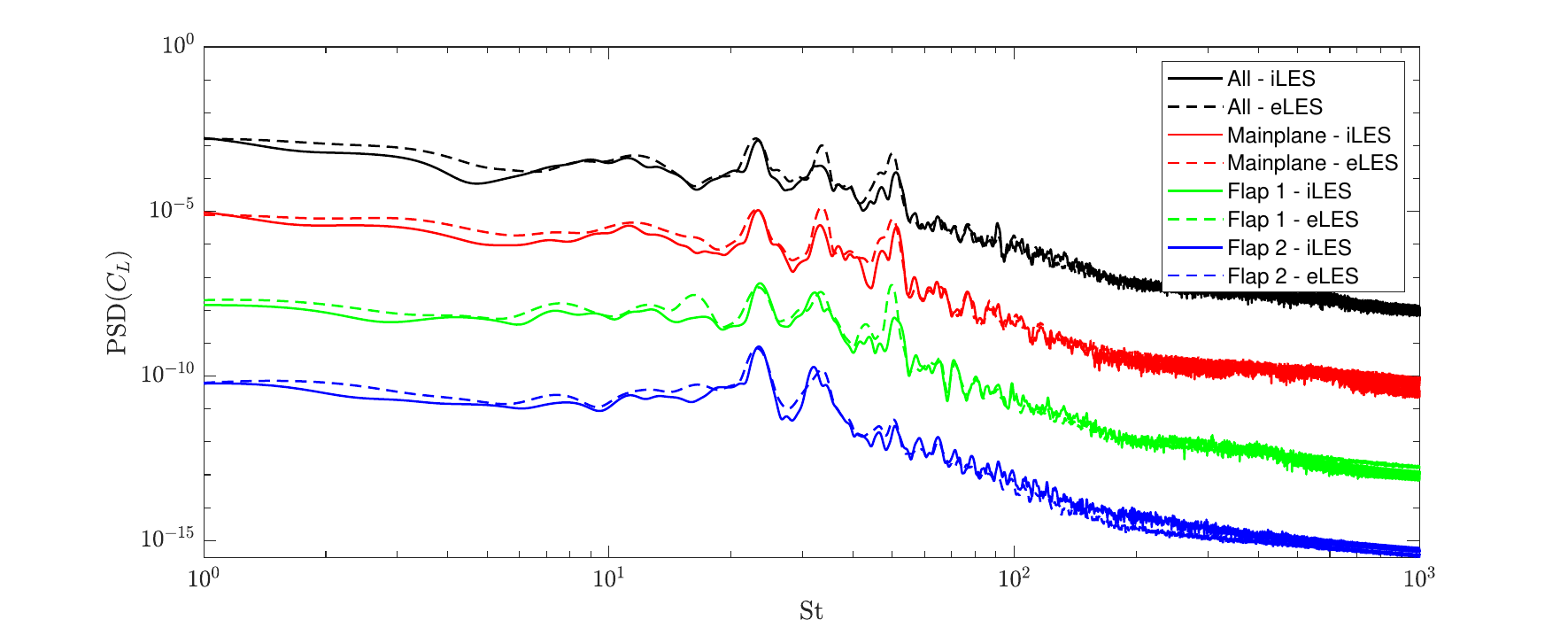}
    \caption{PSD plot of the total lift coefficient and the lift coefficient of each individual element for the iLES and eLES methods of Horses3D.}
   \label{fig:lift_psd}
\end{figure}

As part of the analysis, we are also interested in the predicted values of the force coefficients as they play a pivotal role in the prediction of the performance and the design choices in Formula 1 cars. In Figures \ref{fig:lift_coeff} and \ref{fig:drag_coeff}, we present the time series for the lift and drag coefficients, respectively. After allowing the flow to evolve for 40 CTU we reach a statistically stable solution, as shown in the plots, and we perform an averaging over an additional 8 CTU. A comparison of the averaged quantities is presented in table~\ref{tab:force_coefs}. The difference between the iLES and the eLES results is $\Delta C_{d} = 0.009$ and $\Delta C_{l} = 0.1147$ which is $5\%$ and $1.3\%$ accordingly with respect to the iLES results. In table~\ref{tab:force_coefs} we have also included the results obtained with Nektar++. Although, some discrepancy is expected, the predictions are all within a $6\%$ window from the Nektar++ results. As reported in \cite{slaughter2023large}, for a varying span size such difference in the force coefficients is expected and the results of this work fall within that window. 
\begin{table}[H]
    \centering
    \begin{tabular}{cccc}
         & Horses3D iLES & Horses3D LES & Nektar++ \cite{slaughter2023large}\\ \hline
    $C_{L}$  & -8.6821 & -8.5674 & -8.33  \\
    $C_{D}$  & 0.1725  & 0.1815  & 0.170 \\
    \end{tabular}
    \caption{Average force coefficients for the Horses3D and Nektar++ results.}
    \label{tab:force_coefs}
\end{table}
To further understand the characteristics of the flow and the phenomena that take place we calculate the power spectral density (PSD) of the time series of the $C_{L}$. To do so, we follow \cite{slaughter2023large} and use a signal sampled at each $dt^{*}$. This is calculated using MATLAB's pweltch algorithm with a segment size of 62,500 and an overlapping of $50\%$. In Figure~\ref{fig:lift_psd} we present the PSD for the total lift coefficient $C_{L}^{ALL}$ and for each individual element. 
The PSD of $C_{L}^{ALL}$ shows three peaks for the two simulations. Those are at $St =23,33,50$. These are the dominant frequencies appearing in the PSD of the lift coefficient of the mainplane $C_{L}^{main}$. As proven and described in \cite{slaughter2023large} these are directly related to the LSB mechanism on the suction side. For the first flap $C_{L}^{1^{st}}$, the PSD has the same peaks as those of the mainplane. This is also reported in \cite{slaughter2023large}. Although the scale and the upstream conditions for each element are different, there seems to be a direct relation between the two. For the $2^{nd}$ flap the PSD of $C_{L}^{2^{nd}}$ showcases two peaks at $St =23,33$. The dominant frequencies also appear in the $2^{nd}$ flap in \cite{slaughter2023large}. The trend and the interactions of all three elements presented in \cite{slaughter2023large} seem to be captured in this work. No differences are present between the iLES and the eLES results. However, there are significant differences in the location of the peaks with the dominant modes in \cite{slaughter2023large} being present at $St = 60, 140, 200$. Also, lower mode frequencies at $St =30,40$ are present in \cite{slaughter2023large}. Due to the differences in the setup this discrepancy is expected.

\subsection{Wake Statistics}
\label{sec:wake_stats}
\begin{figure}\centering
\subfloat[ $x/c = 2$]{\label{fig:wake_u_xc2}\includegraphics[width=.49\linewidth]{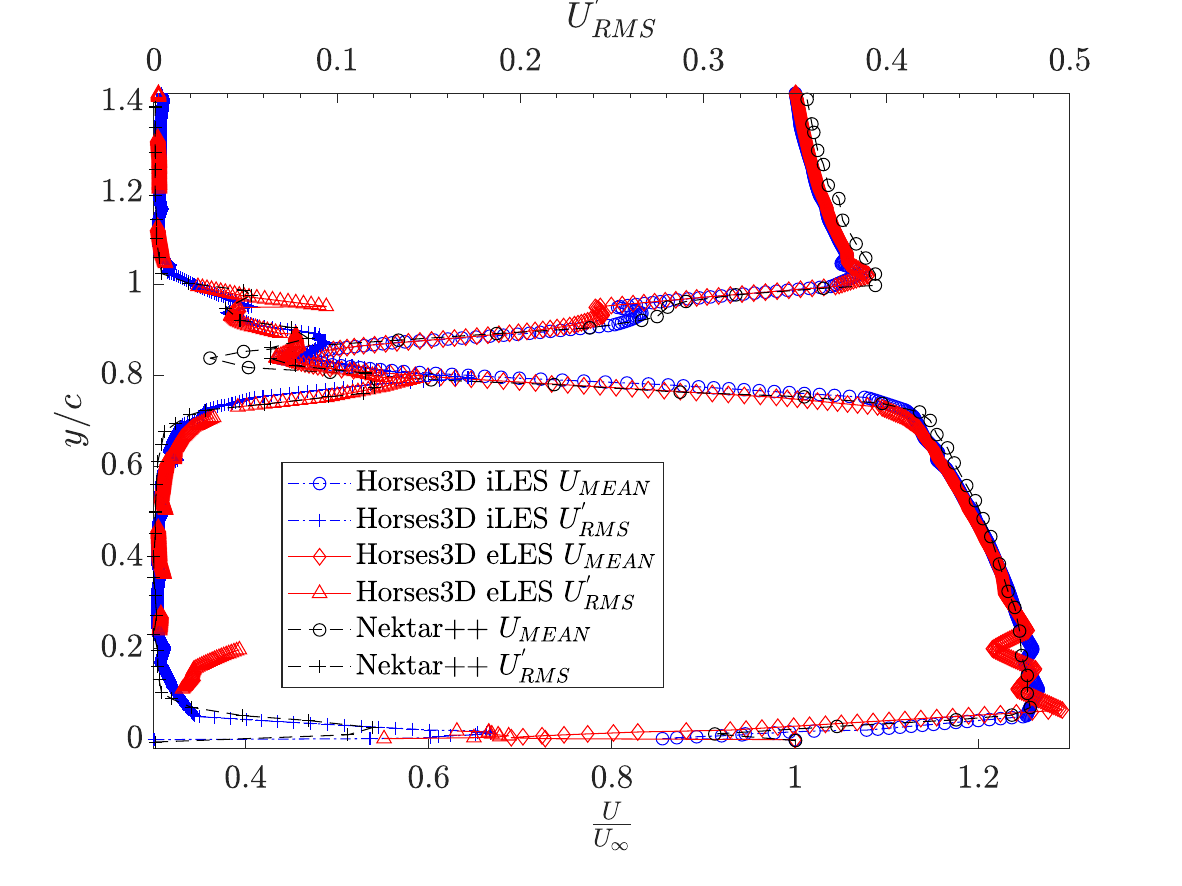}}
\subfloat[$x/c = 2.5$]{\label{fig:wake_u_xc2d5}\includegraphics[width=.49\linewidth]{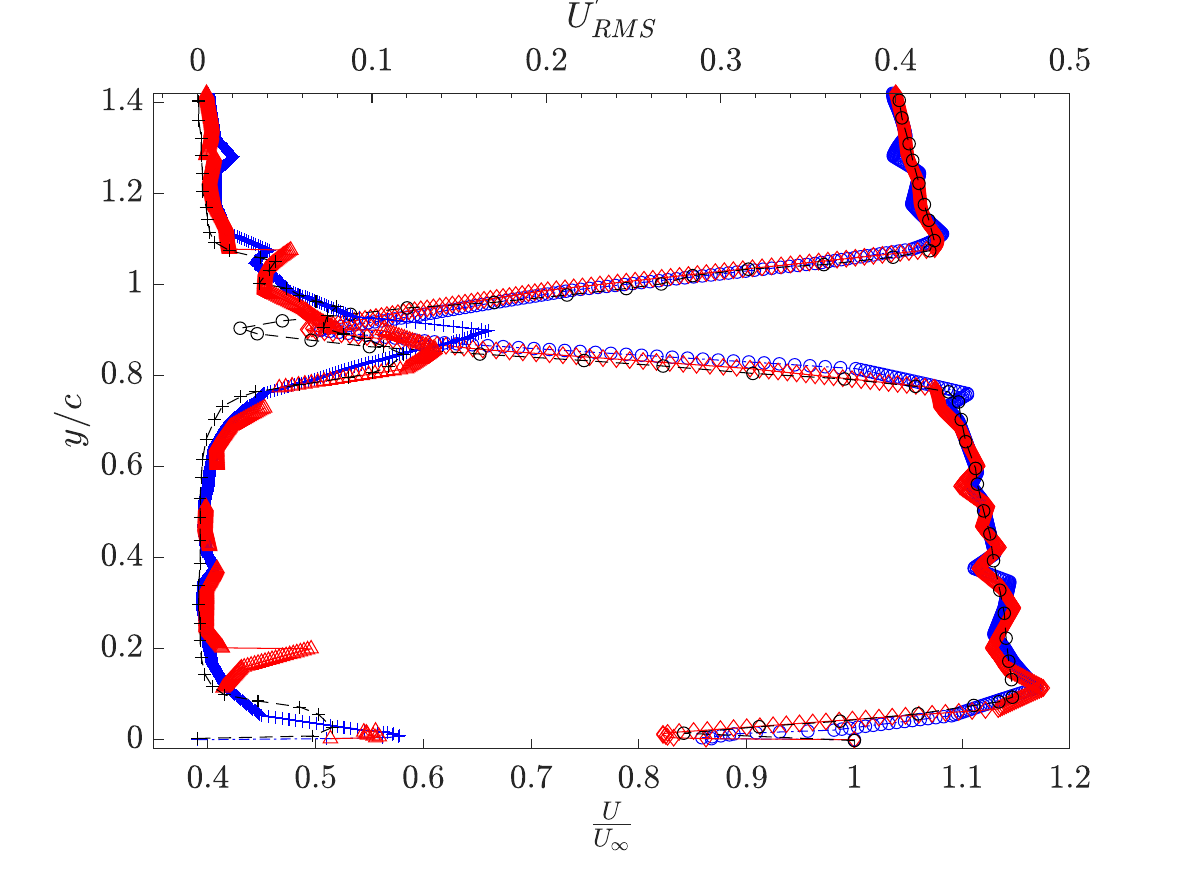}}\par 
\subfloat[$x/c = 4$]{\label{fig:wake_u_xc4}\includegraphics[width=.49\linewidth]{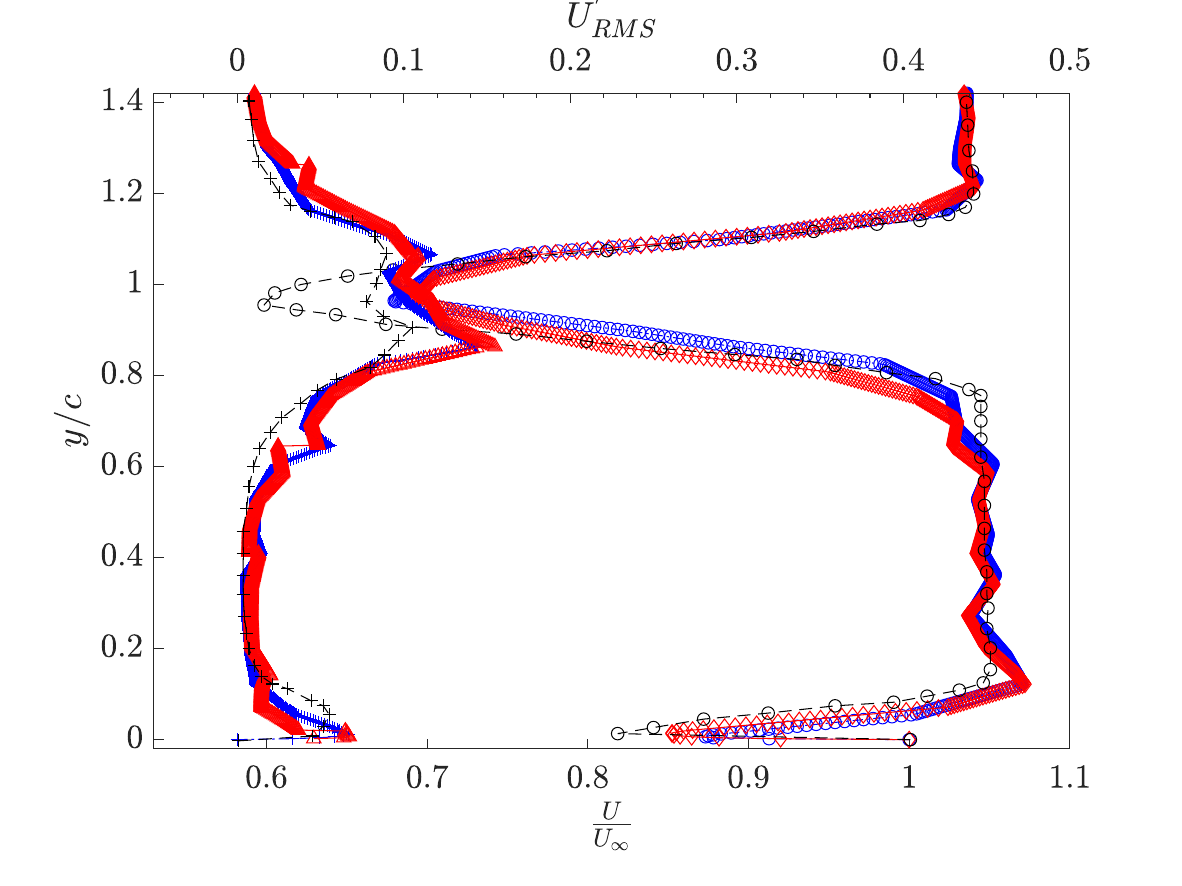}}
\subfloat[$x/c = 6$]{\label{fig:wake_u_xc6}\includegraphics[width=.49\linewidth]{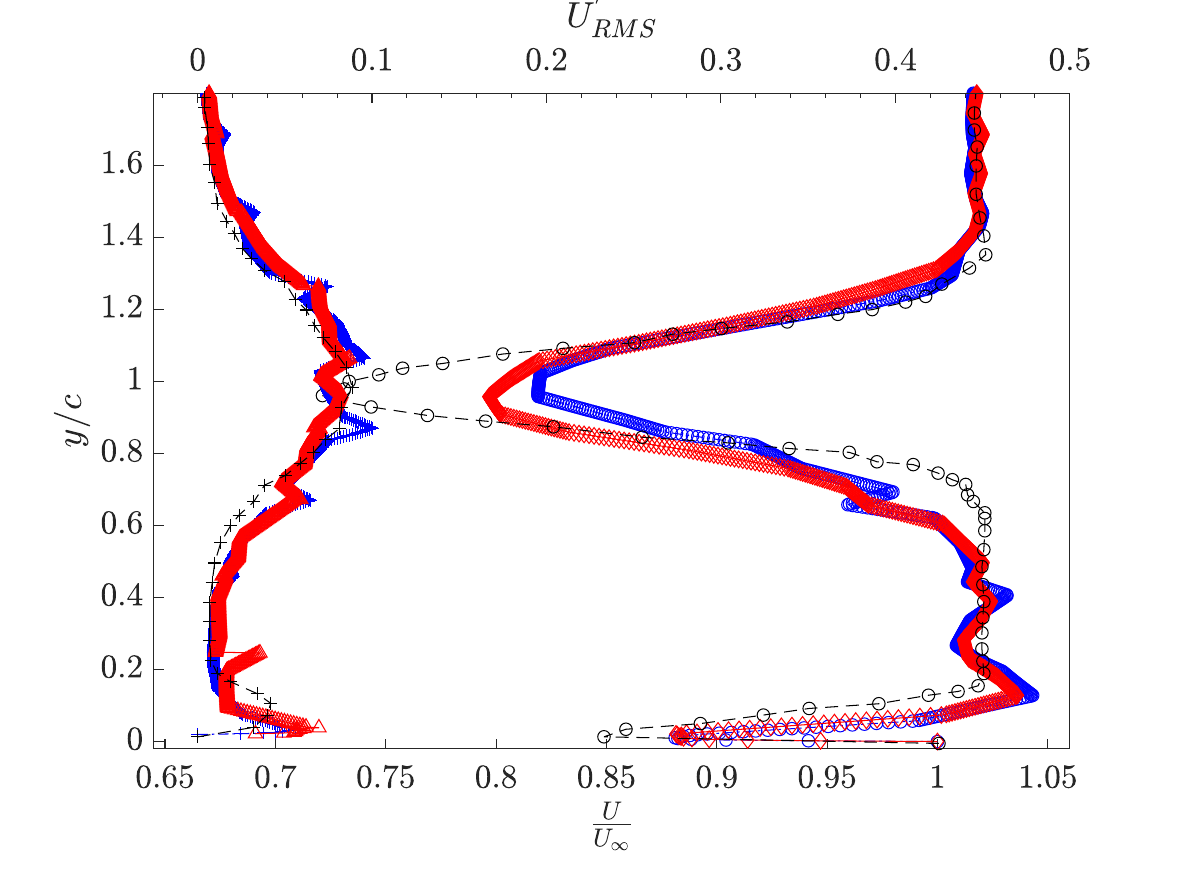}}
    \caption{Comparison of the wake momentum deficit $U/U_{\infty}$ and the $u_{RMS}$ for the three methods at each different station of the wake.}
   \label{fig:wake_u_all}
\end{figure}

In Figure~\ref{fig:wake_u_all}, we present a comparison for the wake momentum deficit and the $u_{RMS}$ at four different stations in the wake. These are $x / c = 2,2.5,4,6$ measured from the leading edge of the mainplane. We present the results of the iLES, eLES simulations from Horses3D and that of Nektar++. In general, for the average momentum deficit $U/U_{\infty}$ and the $u_{RMS}$ match between the iLES and the eLES along the different control stations. For the  $x/c=2$, there is very good match between the data from Netkar++ \cite{slaughter2023large} and the results from Horses3D. The momentum deficit near the floor is also captured properly. A small discrepancy is noted on the point of minimum momentum, as the results from Horses3D appear to have a slightly increased minimum. A good match is also apparent in the $u_{RMS}$ values in this station with the peak value near the floor and within the wake being overestimated. A similar pattern is observed for the $x/c=2.5$ station, with a good overall match. The minimum momentum deficit is overestimated whereas the peak of the $u_{RMS}$ in the wake region is slightly higher.  

As we move further downstream, we observe a more significant deviation between the two datasets. This can be attributed to the coarser mesh used with Horses3D in that region. The exact difference is unknown as this is information is not provided in \cite{slaughter2023large}. For the station $x/c=4$, the momentum deficit near the floor has a difference of $7\%$. A more pronounced difference is also visible in the wake region as the minimum momentum deficit is over-predicted by $8\%$. The $u_{RMS}$ are in agreement, with the peaks in the wake boundaries with the freestream being over-predicted. For the $x/c=6$ station, the difference in the wake momentum deficit is $10\%$. This is again a consequence of the mesh resolution. The $u_{RMS}$ is a good match also in this case. As in \cite{slaughter2023large}, the location of the peak of the momentum deficit appears to be higher for $x/c=4$ compared to that of $x/c=2.5$ and $x/c=6$.

\begin{figure}\centering
\subfloat[ $x/c = 2$]{\label{fig:wake_tke_xc2}\includegraphics[width=.49\linewidth]{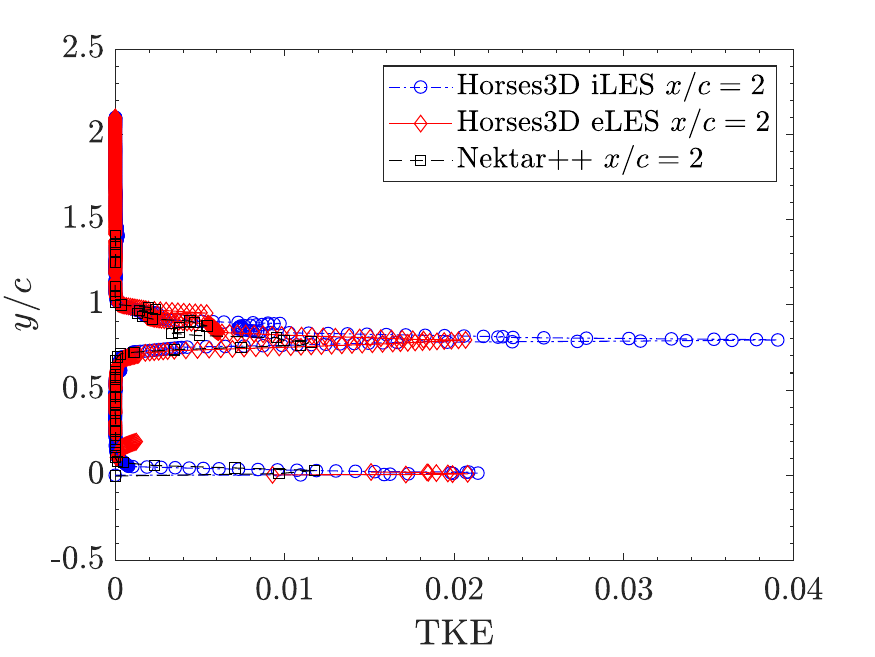}}
\subfloat[$x/c = 2.5$]{\label{fig:wake_tke_xc2d5}\includegraphics[width=.49\linewidth]{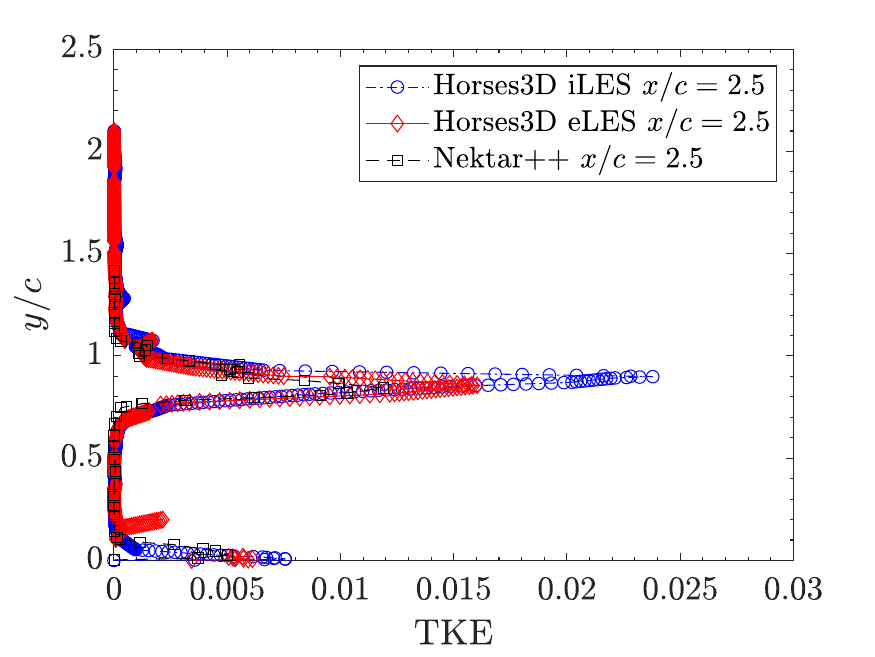}}\par 
\subfloat[$x/c = 4$]{\label{fig:wake_tke_xc4}\includegraphics[width=.49\linewidth]{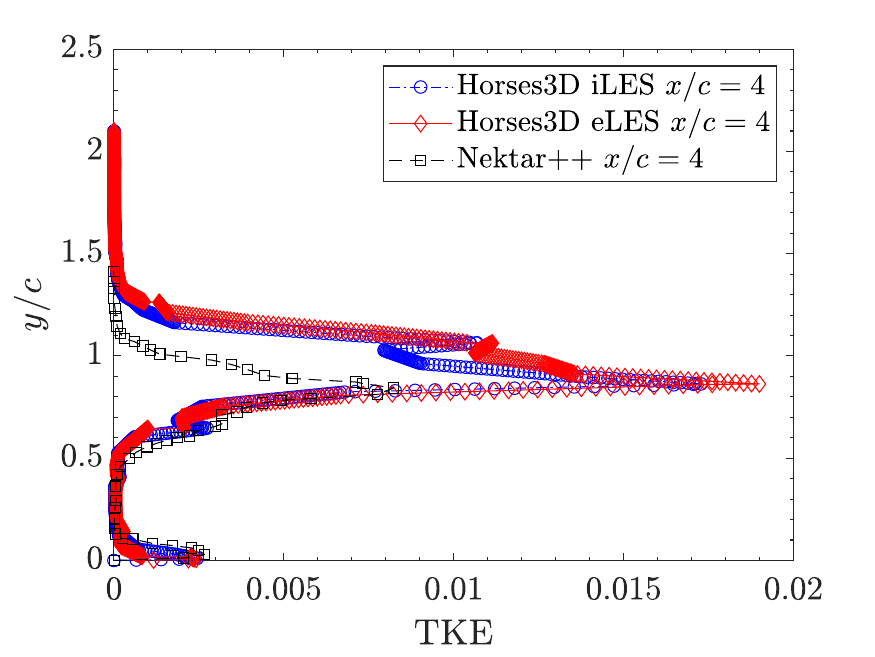}}
\subfloat[$x/c = 6$]{\label{fig:wake_tke_xc6}\includegraphics[width=.49\linewidth]{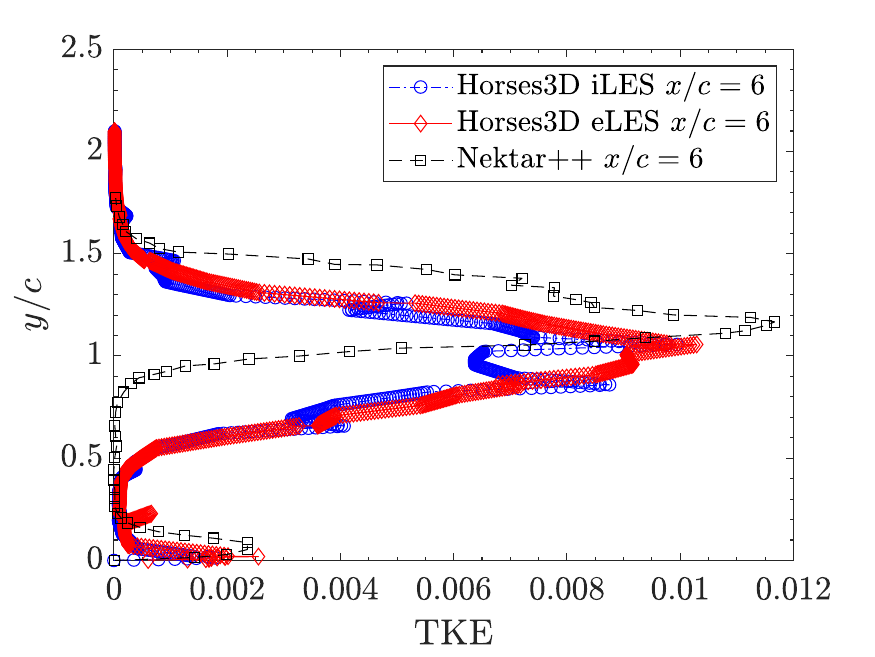}}
    \caption{Comparison of TKE for the three methods at each different station of the wake.}
   \label{fig:wake_tke_all}
\end{figure}

In Figure~\ref{fig:wake_tke_all} we present a comparison of the turbulent kinetic energy (TKE) at each of the four stations. The TKE is defined as $TKE  = \frac{1}{2}\left(\overline{u^{'}u^{'}} + \overline{v^{'}v^{'}} + \overline{w^{'}w^{'}}\right)$. Focusing on Figure~\ref{fig:wake_tke_xc2}, we observed that the peak near the floor and the wake is heavily over-estimated in Horses3D. The exact source of this is unknown and it can be attributed to the coarser resolution of the mesh and the lower spanwise length. However, within the wake region we can identify three peaks which correspond to the wakes of each individual element. The peaks of the mainplane and the $1^{st}$ flap are over-estimated, while there is a very good agreement for the peak of the $2^{nd}$ flap. The location of the peaks matches that reported in \cite{slaughter2023large}. The iLES results significantly overestimate the TKE compared to the eLES and the Nektar++ results. A similar trend is observed in Figure~\ref{fig:wake_tke_xc2d5} for the  $x/c=2.5$ station. Although in this case there is a good match for the peak of the $1^{st}$ and the $2^{nd}$ flaps.

In Figure~\ref{fig:wake_tke_xc4}, as we examine a region further away from the front wing section, the iLES and eLES results seem to come to an agreement, with the differences being less pronounced compared to the stations closer to the wing. However, we observe more significant differences between the results of Horses3D and Nektar++ in the region of the wake. In the results of this work, two peaks are present in the TKE. The presence of the second peak is probably due to insufficient mixing, caused by the lack of resolution in that region. Also, in this case, the peak of TKE is over-estimated. Lastly, for the station $x / c = 6$, the results in Figure~\ref{fig:wake_tke_xc6} show the highest discrepancy, as expected. In contrast to all the previous results, the TKE is under-predicted and its the only station for which the location of the peaks within the wake is not in agreement with \cite{slaughter2023large}. 

\subsection{Computational cost}

Finally, we conduct a comparative analysis of the explicit (Vreman) and implicit (KG) Large Eddy Simulation formulations within the Horses3D framework. 
Our primary focus is to determine the maximum stable time step achievable using the explicit time marching RK3 scheme for each LES formulation and in addition assess the computational costs.

To ascertain the maximum time step, we initiated computations from a uniform initial condition and progressively increased the CFL number from 0.5 in increments of 0.1. Running 100 iterations for both formulations revealed that the explicit LES model (Vreman) experiences instability for CFLs above $CFL_{Vreman-limit}=0.6$, 
with an associated time step of $dt=3.25\times 10^{-6}$. On the contrary, the KG model remained stable for CFLs up to  $CFL_{KG-limit}=1.0$ and associated $dt=4.94\times 10^{-6}$. 
Consequently, we observed that the KG model allows for a time step approximately 50\% greater than that of Vreman. The increased stability of the KG scheme is not associated with the turbulent closure but to the use of Gauss or Gauss-Lobatto points for integration, the latter allowing larger CFLs. These findings align with those presented in \cite{doi:10.1137/100807211}, indicating a time step capability up to 2 times greater in Gauss-Lobatto compared to Gauss. 

Regarding the computational cost, we computed the average cost for 100 iterations (and for three separate runs) on an Intel(R) Xeon(R) Gold 6248R CPU @ 3.00GHz processor, equipped with 2 sockets and 24 cores per socket. The simulations were conducted using MPI parallelization with 48 processes. Our findings are reported in table \ref{tab:my-table} and show that while both formulations exhibit similar computational cost, KG exhibits marginally lower cost per iteration when compared to Vreman, on this specific computing system. The cost per iteration of KG and Vreman is 10\% higher than a vanilla DGSEM solver (Gauss nodes, without splitting -no energy/entropy stable formulation- and without explicit turbulence model).

Additionally, we report the computational cost per CTU, which comprises both effects (cost per iteration and maximum stable time step). We can see in table \ref{tab:my-table} that if the physical problem permits using the maximum stable time step, implicit LES is around 50\% faster in Horses3D for this Formula 1 front wing using a polynomial order 4 (5th order accurate). 
Note that the wall time the simulation can be decreased by increasing the number of nodes used. Horses3D has been tested and presents good scalability properties up to $10,000$ nodes as presented in Appendix \ref{app:scalability}. 
\begin{table}[H]
    \centering
    \begin{tabular}{ccccc}
         & \textbf{CFL} & 
         \textbf{\begin{tabular}[c]{@{}l@{}}max time step\\ (non dimensional)\end{tabular}} & 
         \textbf{\begin{tabular}[c]{@{}l@{}} CPU cost (s) \\ per iteration per node \end{tabular}} & 
         \textbf{\begin{tabular}[c]{@{}l@{}}CPU cost (h) \\ per CTU per node\end{tabular}} \\ \hline                                             
    \textbf{eLES}  & $0.6$ & $3.25\times 10^{-6}$ & $1.42$ &   $121.1$ \\ 
    \textbf{iLES}  & $1.0$ & $4.94\times 10^{-6}$ & $1.38$ &   $77.5$   \\ 
    \end{tabular}
    \caption{Maximum stable time step and cost per iteration and per CTU for explicit (Vreman) LES and implicit (KG) LES.}
    \label{tab:my-table}
\end{table}

\section{Conclusions}
We have performed explicit and implicit simulations and analysis of a simplified configuration of a real Formula 1 front wing under nominal conditions using a high-order discontinuous Galerkin method. The simulations have been carried out using the Horses3D framework \cite{HORSES3D} using a polynomial order 4 ($5th$ order accuracy). We have presented a direct qualitative and quantitative comparison with the results obtained using Nektar++ \cite{slaughter2023large}. In this study, the number of DoF is approximately one order of magnitude lower than that used in the reference. Both numerical techniques implicit and explicit LES provide accurate results and compare favourably to the reference data. 

With respect to the accuracy of the prediction, the deviation for the time-averaged force coefficients compared to \cite{slaughter2023large} is below 5\% in all cases. The results near and on the surface of each element are also in very good agreement. The $C_{p}$ and $C_{f}$ distributions across each element are almost identical between the data sets, with only some minor differences. The existence and position of each LSB is accurately predicted with the current methods, with a slightly larger bubble for the explicit LES simulation. 
Furthermore, in the wake region, where we have under-resolved conditions we observe that the wake statistics are in good agreement, given the number of DoF of the simulations in this work. The effects of the under-resolution are more evident in the calculations of the TKE in the different stations of the wake. 

Regarding the explicit and implicit LES comparison, we conclude that:
\begin{itemize}
    \item The implicit LES shows to better capture transition and allows for larger time steps at a similar cost per iteration, and conclude that this formulation is very attractive for complex physics simulations, as exemplified by the Formula 1 front wing.
    \item The implicit LES scheme is robust, throughout the tests, even in the absence of stabilising techniques such as sub-grid scale modelling. This is a highly desired behavior towards the effort for the industrialisation of said methods. 
\end{itemize}


The methods presented in this paper are general and can be used to arbitrary 3D geometries and meshes without further modification. This work constitutes an additional step towards the better understanding of the applicability of high-order methods on industrial automotive test cases and an insight on how they could be used efficiently depending on the desired targets of a CFD study. 

\appendix

\section{Horses3D strong scalability}
\label{app:scalability}

We have performed scalability tests in LUMI-C thanks to the EHPC-BEN program. Two test cases are examined to assess strong scalability of Horses3D. The first case focuses on the Taylor-Green Vortex (TGV) \cite{TaylorGreen} problem at a Reynolds number of 1600. The TGV setup involves simulating laminar to turbulent flow in a three-dimensional box, utilizing varying h-meshes and polynomial orders. The second case involves simulating airflow around a NACA0012 airfoil using a 56088-element h-mesh with polynomial orders P = 4 and P = 6. Unlike the TGV problem, this scenario includes complex boundary conditions, such as inflow, outflow, periodic, and no-slip wall conditions. The purpose is to test the boundary condition impact on the scalability of Horses3D. A summary of the test cases conducted is presented in table \ref{tab:scalability}.

\begin{table}[h]
\begin{tabular}{llllll}
       & TEST     & PARALELIZATION STRATEGY & ELEMENTS & P &  \\
TEST 1 & TGV      & MPI ONLY                & $128^3$  & 3 &  \\
TEST 2 & TGV      & MPI ONLY                & $64^3$   & 7 &  \\
TEST 3 & NACA0012 & OpenMP + MPI            & 56088    & 4 &  \\
TEST 4 & NACA0012 & OpenMPI + MPI            & 56088    & 6 & 
\end{tabular}
\caption{Summary of scalability test cases}
\label{tab:scalability}
\end{table}

\begin{figure}[h]\centering
\includegraphics[width=.75\linewidth]{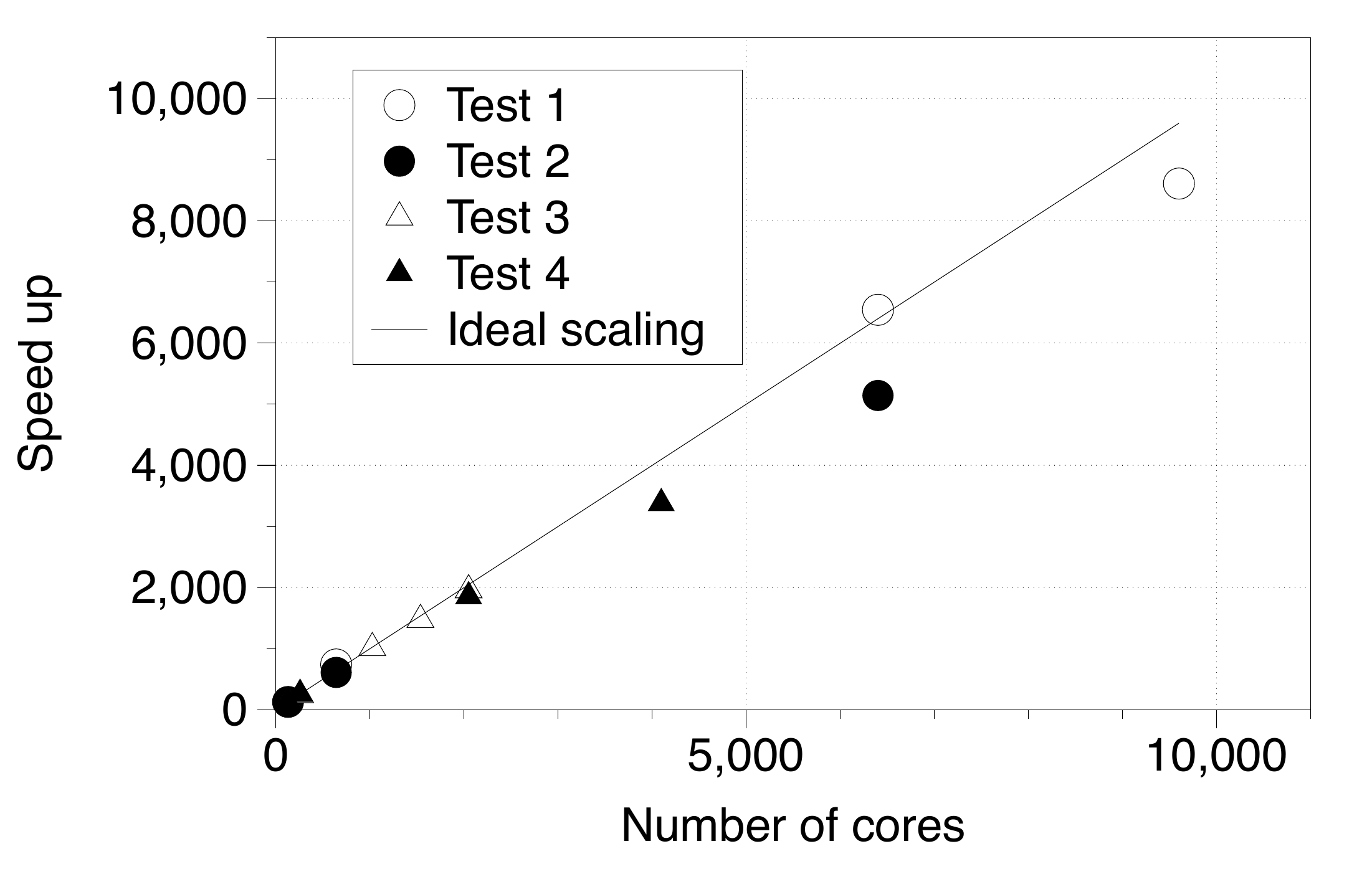}
    \caption{Strong scalability tests}
   \label{fig:scalability}
\end{figure}

\subsection{Taylor Green Vortex (TGV): MPI strong scaling}
We perform strong scalability simulations using only MPI processors (up to 9000). We use a h-mesh composed of $64^3$ and $128^3$ elements that we enrich with polynomials P = 3 and P = 7. Total number of DoF is computed as $Elements\times(P+1)^3$. Figure \ref{fig:scalability} shows that very good scaling (almost ideal) is obtained when the loading is high (mesh with $128^3$ - P=3) for up to 9000 processors and some minor degradation if a small mesh is used. Similar performance is obtained with the $64^3$ mesh when varying the polynomial order. This last remark shows the potential of the solver since much more resolution is obtained when using higher polynomials whilst the scalability remains similar. 

\subsection{Flow around a NACA0012 airfoil – mixed OpenMP – MPI strong scaling}
Horses3D also permits hybrid OpenMP + MPI paralelization. Based on preliminary simulations, we find that the sweet spot for OpenMP is with 4-8 OpenMP threads but experiences significant degradation beyond this range. In all instances, we utilize the settings "export OMP\_PROC\_BIND=close" and "export OMP\_PLACES=cores" to optimize performance.
We check the MPI performance when fixing to 8 the OpenMP threads. The findings are succinctly presented in Figure \ref{fig:scalability}. It is evident that nearly perfect scaling is achieved with 2,000 threads, although a slight degradation in performance becomes apparent when scaling up to 4,000 threads. This last case shows that including the typical boundary conditions and LES turbulent models do not damage the scaling of Horses3D. Furthermore, we observe excellent scaling (almost ideal) when combining moderate number of OpenMP (up to 8 threads) with large number of MPI processors.

\section{Compressible Navier-Stokes}
\label{sec:cNS}
In this work we solve the 3D Navier-Stokes for iLES simulations and we supplement the equations with the Vreman LES model for eLES simulations. The 3D Navier-Stokes when including the  Vreman model can be compactly written as:
\begin{equation}
\boldsymbol{u}_t+\nabla\cdot\ssvec{{F}}_e = \nabla\cdot\ssvec{F}_{v,turb},
\label{eq:compressibleNScompact}
\end{equation}
where $\boldsymbol{u}$ is the state vector of large scale resolved conservative variables $\boldsymbol{u} = [ \rho , \rho v_1 , \rho v_2 , \rho v_3 , \rho e]^T$, $\ssvec{F}_e$ are the inviscid, or Euler fluxes,
\begin{equation}
\ssvec{F}_e = \left[\begin{array}{ccc} \rho v_1 & \rho v_2 & \rho u_3 \\
                                                                                \rho v_1^2 + p & \rho v_1v_2 & \rho v_1v_3 \\
                                                                                	\rho v_1v_2 & \rho v_2^2 + p & \rho v_2v_3 \\
                                                                                	\rho v_1v_3 & \rho v_2v_3 & \rho v_3^2 + p \\
                                                                                	\rho v_1 H & \rho v_2 H & \rho v_3 H
\end{array}\right],
\end{equation}
where $\rho$, $e$, $H=E+p/\rho$, and $p$ are the large scale density, total energy, total enthalpy and pressure, respectively, and $\vec{v}=[v_1,v_2,v_3]^T$ is the large scale resolved velocity components. Additionally, $\ssvec{F}_{v,turb}$ and defines the viscous and turbulent fluxes,
\begin{equation}
\ssvec{F}_{v,turb}= \left[\begin{array}{ccc}0 & 0 & 0\\
\tau_{xx} & \tau_{xy} & \tau_{xz} \\
\tau_{yx} & \tau_{yy} & \tau_{yz} \\
\tau_{zx} & \tau_{zy} & \tau_{zz} \\
\sum_{j=1}^3 v_j\tau_{1j} + \kappa T_x& \sum_{j=1}^3 v_j\tau_{2j} + \kappa T_y& \sum_{j=1}^3 v_j\tau_{3j} + \kappa T_z
\end{array}\right],
\label{eq:viscousfluxes}
\end{equation}
where $\kappa$ is the thermal conductivity, $T_x, T_y$ and $T_z$ denote the temperature gradients and the stress tensor $\boldsymbol{\tau}$ is defined as $\boldsymbol{\tau} = (\mu+\mu_t)(\nabla \vec{v} + (\nabla \vec{v})^T) - 2/3(\mu+\mu_t) \boldsymbol{I}\nabla\cdot\vec{v}$, with $\mu$ the dynamic viscosity, $\mu_t$ the turbulent viscosity (in this work defined through the Vreman 
model) and $\boldsymbol{I}$ the three-dimensional identity matrix. 
%
The dynamic turbulent viscosity using the Vreman \cite{vreman2004eddy} model is given by: 
\begin{equation}
\begin{split}
&\mu_t = C_v \rho\sqrt{\frac{B_\beta}{\alpha_{ij}\alpha_{ij}}},\\
&\alpha_{ij} = \frac{\partial{u}_j}{\partial{x}_i},\\
&\beta_{ij} = \Delta^2\alpha_{mi}\alpha_{mj},\\
&\Delta=\frac{V^{1/3}}{N+1},\\
&B_\beta = \beta_{11}\beta_{22} -\beta_{12}^2 +\beta_{11}\beta_{33} -\beta_{13}^2 +\beta_{22}\beta_{33} -\beta_{23}^2,
\end{split}
\label{eq-iLES:LES_vreman}
\end{equation}
 
\noindent where $C_v=0.07$ is the constant of the model, $V$ is the volume of the element and $N$ is the polynomial order of the approximation. The Vreman LES model adjusts the model parameters based on the local flow characteristics and automatically reduces the turbulent viscosity in laminar, transitional, and near-wall regions allowing to capture the correct physics.
\newpage

\section*{Declaration of Competing Interest}
The authors declare that they have no known competing financial interests or personal relationships that could have appeared to influence the work reported in this paper.
\section*{Acknowledgments}
Gonzalo Rubio and Esteban Ferrer  acknowledge the funding received by the Grant DeepCFD (Project No. PID2022-137899OB-I00) funded by MCIN\slash AEI\slash 10.13039\slash 501100011033/ and by ERDF A way of making Europe. 
Esteban Ferrer aknowledges the Comunidad de Madrid and Universidad Politécnica de Madrid for the Young Investigators award: APOYO-JOVENES-21-53NYUB-19-RRX1A0. 
Gerasimos Ntoukas and Esteban Ferrer would like to thank the European Union’s Horizon 2020 Research and Innovation Program under the Marie Skłodowska-Curie grant agreement No 813605 for the ASIMIA ITN-EID project. 
The author thankfully acknowledges the computer resources at MareNostrum and the technical support provided by Barcelona Supercomputing Center (RES-IM-2022-3-0023). 
The authors acknowledge the EuroHPC Joint Undertaking for awarding this project access to the EuroHPC supercomputer LUMI, hosted by CSC (Finland) and the LUMI consortium through a EuroHPC Benchmark Access call (EHPC-BEN-2023B05-015).
The authors gratefully acknowledge the Universidad Politécnica de Madrid (www.upm.es) for providing computing resources on Magerit Supercomputer.


\bibliographystyle{vancouver}
\bibliography{biblio }
\end{document}